\documentclass[aos,MSNbibl,nameyear,seceqn,dvips]{arximspdf}
\usepackage{mathbh,dcolumn}
\usepackage{graphicx}

\doi{10.1214/12-AOS1049} 
\volume{40}
\issue{5}
\pubyear{2012}
\firstpage{2733}
\lastpage{2763}

\makeatletter
\newcolumntype{d}[1]{D{.}{.}{#1}}
\renewcommand{\citep}[1]{(\citeauthor{#1} \citeyear{#1})}
\newtheorem{thmm}{Theorem}
\newtheorem{cor}[thmm]{Corollary}
\newcommand{\argmax}{\operatorname{argmax}}
\makeatother

\begin{document}
\begin{frontmatter}

\title{Optimal weighted nearest neighbour classifiers\thanksref{T1}}
\thankstext{T1}{Supported by the Leverhulme Research Fellowship and an
EPSRC Early Career Fellowship.}
\runtitle{Weighted nearest neighbour classifiers}

\begin{aug}
\author{\fnms{Richard J.} \snm{Samworth}\corref{}\ead[label=e1]{r.samworth@statslab.cam.ac.uk}\ead[label=u1,url]{http://www.statslab.cam.ac.uk/\textasciitilde rjs57}}
\runauthor{R.~J. Samworth}
\affiliation{University of Cambridge}
\address[A]{Statistical Laboratory\\
University of Cambridge\\
Wilberforce Road\\
Cambridge \\
CB3 0WB \\
United Kingdom\\
\printead{e1}\\
\printead{u1}} 
\end{aug}

\received{\smonth{6} \syear{2012}}
\revised{\smonth{8} \syear{2012}}

%
\begin{abstract}
We derive an asymptotic expansion for the excess risk (regret) of a
weighted nearest-neighbour classifier. This allows us to find the
asymptotically optimal vector of nonnegative weights, which has a
rather simple form. We show that the ratio of the regret of this
classifier to that of an unweighted $k$-nearest neighbour classifier
depends asymptotically only on the dimension~$d$ of the feature
vectors, and not on the underlying populations. The improvement is
greatest when $d=4$, but thereafter decreases as $d \rightarrow\infty
$. The popular bagged nearest neighbour classifier can also be regarded
as a weighted nearest neighbour classifier, and we show that its
corresponding weights are somewhat suboptimal when $d$ is small (in
particular, worse than those of the unweighted $k$-nearest neighbour
classifier when $d=1$), but are close to optimal when $d$ is large.
Finally, we argue that improvements in the rate of convergence are
possible under stronger smoothness assumptions, provided we allow
negative weights. Our findings are supported by an empirical
performance comparison on both simulated and real data sets.
\end{abstract}

%
\begin{keyword}[class=AMS]
\kwd{62G20}
\end{keyword}

\begin{keyword}
\kwd{Bagging}
\kwd{classification}
\kwd{nearest neighbours}
\kwd{weighted nearest neighbour classifiers}
\end{keyword}

\end{frontmatter}
%
\section{Introduction}
\label{SecIntro}

Supervised classification, also known as pattern recognition, is a
fundamental problem in Statistics, as it represents an abstraction of
the decision-making problem faced by many applied practitioners.
Examples include a doctor making a medical diagnosis, a handwriting
expert performing an authorship analysis, or an email filter deciding
whether or not a message is genuine.

Classifiers based on nearest neighbours are perhaps the simplest and
most intuitively appealing of all nonparametric classifiers. The
$k$-nearest neighbour classifier was originally studied in the seminal
works of \citet{FixHodges1951} [later republished as \citet
{FixHodges1989}] and \citet{CoverHart1967}, but it retains its
popularity today. Surprisingly, it is only recently that detailed
understanding of the nature of the error probabilities has emerged
[\citet{HPS2008}].

Arguably the most obvious defect with the $k$-nearest neighbour
classifier is that it places equal weight on the class labels of each
of the $k$ nearest neighbours to the point $x$ being classified.
Intuitively, one would expect improvements in terms of the
misclassification rate to be possible by putting decreasing weights on
the class labels of the successively more distant neighbours.

The first purpose of this paper is to describe the asymptotic structure
of the difference between the misclassification rate (risk) of a
weighted nearest neighbour classifier and that of the optimal Bayes
classifier for classification problems with feature vectors in $\mathbb
{R}^d$. Theorem~\ref{ThmMain} in Section~\ref{SecMain} below shows
that, subject to certain regularity conditions on the underlying
distributions of each class and the weights, this excess risk (or \emph
{regret}) asymptotically decomposes as a sum of two dominant terms, one
representing bias and the other representing variance. For simplicity
of exposition, we will deal initially with binary classification
problems, though we also indicate the appropriate extension to general
multicategory problems.

\begin{figure}

\includegraphics{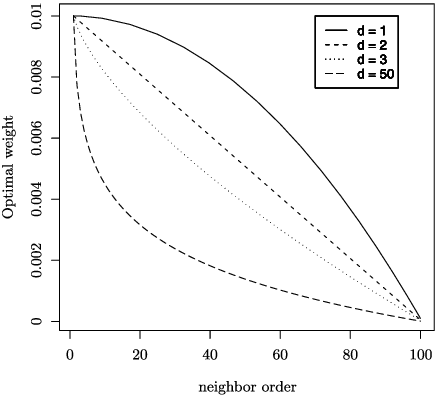}

\caption{Optimal weight profiles at different dimensions.
Here, $k^*=100$, and the figure displays the positive weights
in~(\protect\ref{EqOptWeights}),
scaled to have the same weight on the nearest neighbour at each
dimension.}\label{FigDecWeights}
\end{figure}

Our second contribution, following on from the first, is to derive the
vector of nonnegative weights that is asymptotically optimal in the
sense of minimising the misclassification rate; cf. Theorem~\ref
{ThmOptWeights}. In fact this asymptotically optimal weight vector has
a relatively simple form: let $n$ denote the sample size and let
$w_{ni}$ denote the weight assigned to the $i$th nearest neighbour
(normalised so that $\sum_{i=1}^n w_{ni} = 1)$. Then the optimal choice
is to set $k^* = \lfloor B^* n^{4/(d+4)} \rfloor$ [an explicit
expression for $B^*$ is given in~(\ref{Eqkstar}) below] and then let
%
\begin{equation}
\label{EqOptWeights} w_{ni}^* = \cases{\displaystyle
\frac{1}{k^*} \biggl[1 + \frac{d}{2} - \frac
{d}{2(k^*)^{2/d}}\bigl
\{i^{1+2/d} - (i-1)^{1+2/d}\bigr\} \biggr], \vspace*{2pt}\cr
\quad\hspace*{20pt} $\mbox{for $i=1,
\ldots,k^*$},$
\vspace*{2pt}\cr
0, \qquad $\mbox{for $i=k^*+1,\ldots,n$.}$}
\end{equation}
Thus, in the asymptotically optimal weighting scheme, only a proportion
$O(n^{-d/(d+4)})$ of the weights are positive. The maximal weight is
almost $(1+d/2)$ times the average positive weight, and the discrete
distribution on $\{1,\ldots,n\}$ defined by the asymptotically optimal
weights decreases in a concave fashion when $d=1$, in a linear fashion
when $d=2$ and in a convex fashion when $d \geq3$; see Figure~\ref
{FigDecWeights}. When $d$ is large, about $1/e$ of the weights are
above the average positive weight.

Another consequence of Theorem~\ref{ThmOptWeights} is that $k^*$ is
bigger by a factor of $ \{\frac{2(d+4)}{d+2} \}^{d/(d+4)}$ than
the asymptotically optimal choice of $k$ for traditional, unweighted
$k$-nearest neighbour classification. It is notable that this factor,
which is around 1.27 when $d=1$ and increases towards 2 for large $d$,
does not depend on the underlying populations. This means that there is
a natural correspondence between any unweighted $k$-nearest neighbour
classifier and one of optimally weighted form, obtained by multiplying
$k$ by this dimension-dependent factor to obtain the number $k'$ of
positive weights for the weighted classifier, and then using the
weights given in~(\ref{EqOptWeights}) with $k'$ replacing $k^*$.

\begin{figure}

\includegraphics{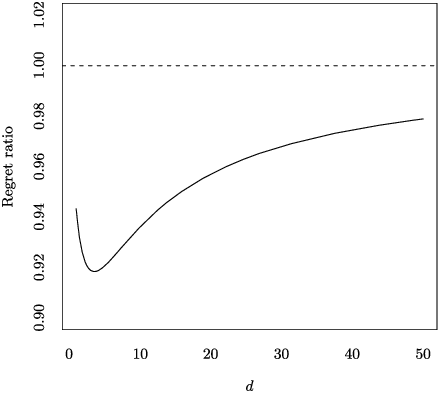}

\caption{Asymptotic ratio of the regret of the optimally weighted
nearest neighbour classifier to that of the optimal $k$-nearest neighbour
classifier,
as a function of the dimension $d$ of the feature vectors.}
\label{FigAsympImp}
\end{figure}

In Corollary~\ref{CorRegretRatio} we describe the asymptotic improvement in the excess
risk that is attainable using the procedure described in the previous
paragraph. Since the rate of convergence to zero of the excess risk is
$O(n^{-4/(d+4)})$ in both cases, the improvement is in the leading
constant, and again it is notable that the asymptotic improvement does
not depend on the underlying populations. The improvement is relatively
modest, which goes some way to explaining the continued popularity of
the (unweighted) $k$-nearest neighbour classifier. Nevertheless, for $d
\leq15$, the improvement in regret is at least 5\%, though it is
negligible as $d \rightarrow\infty$; the greatest improvement occurs
when $d=4$, and here it is just over 8\%. See Figure~\ref{FigAsympImp}.

Another popular way of improving the performance of a classifier is by
\emph{bagging} [Breiman (\citeyear{Breiman1996,Breiman1999})]. Short for ``bootstrap
aggregating'', bagging involves combining the results of many
empirically simulated predictions. Empirical analyses [e.g., \citet
{Steele2009}], have reported that bagging can result in improvements
over unweighted $k$-nearest neighbour classification. Moreover, as
explained by \citet{BCG2010}, understanding the properties of the
bagged nearest neighbour classifier is also of interest because they
provide insight into \emph{random forests} [\citet{Breiman2001}]. Random
forest algorithms have been some of the most successful ensemble
methods for regression and classification problems, but their
theoretical properties remain relatively poorly understood. When
bagging the nearest neighbour classifier, we can draw resamples from the
data either with- or without-replacement. We treat the ``infinite
simulation'' case, where both versions take the form of a weighted
nearest neighbour classifier with weights decaying approximately
exponentially on successively more distant observations from the point
being classified [\citet{HallSamworth2005}, \citet{BCG2010}]. The crucial choice
is that of the resample size, or equivalently the sampling fraction,
that is, the ratio of the resample size to the original sample size. In
Section~\ref{SecBagged}, we describe the asymptotically optimal
resample fraction (showing in particular that it is the same for both
with- and without-replacement sampling) and compare its regret with
those of the weighted and unweighted $k$-nearest neighbour classifiers.

In Section~\ref{SecNegWeights}, we consider the problem of choosing
optimal weights without the restriction that they should be
nonnegative. The situation here is somewhat analogous to the use of
higher order kernels for classifiers based on kernel density estimates
of each of the population densities. In particular, subject to
additional smoothness assumptions on the population densities, we find
that powers of $n$ arbitrarily close to the ``parametric rate'' of
$O(n^{-1})$ for the excess risk are attainable. Section~\ref
{SecEmpirical} presents the results of an empirical performance
comparison of different classifiers studied in the paper, and shows
that the asymptotic theory predicts the empirical performance well. The
main steps in the proof of Theorem~\ref{ThmMain} are given in the
\hyperref[app]{Appendix}; the remaining details can be found in the supplementary
material [\citet{Samworth2012}], along with the other proofs and some
ancillary material.

Classification has been the subject of several book-length treatments,
including \citet{Hand1981}, \citet{DGL1996} and \citet{Gordon1999}. In
particular, classifiers based on nearest neighbours form a central theme
of \citet{DGL1996}. The review paper by \citet{BBL2005} contains 243
references and provides a thorough survey of the classification
literature up to 2005. More recently, \citet{AudibertTsybakov2007} have
discussed the relative merits of \emph{plug-in} classifiers (a family
to which weighted nearest neighbour classifiers belong) and classifiers
based on \emph{empirical risk minimisation}, such as support vector
machines [\citet{CortesVapnik1995}, \citet{BBM2008}, \citet{SteinwartChristmann2008}].

Weighted nearest neighbour classifiers were first studied by \citet
{Royall1966}; see also \citet{BaileyJain1978}. \citet{Stone1977} proved
that if $\max_{1 \leq i \leq n} w_{ni} \rightarrow0$ as $n \rightarrow
\infty$ and $\sum_{i=1}^k w_{ni} \rightarrow1$ for some $k = k_n$ with
$k/n \rightarrow0$ as $n \rightarrow\infty$, then risk of the
weighted nearest neighbour classifier converges to the risk of the Bayes
classifier; see also \citet{DGL1996}, page 179. As mentioned above,
this work attempts to study the difference between these risks more
closely. Weighted nearest neighbour classifiers are also related to
classifiers based on kernel estimates of each of the class densities;
see, for example, the review by \citet{RaudysYoung2004}, as well as \citet
{HallKang2005}. The $O(n^{-4/(d+4)})$ rates of convergence obtained in
this paper for nonnegative weights are the same as those obtained by
\citet{HallKang2005} under similar twice-differentiable conditions with
second-order kernel estimators of the class densities. %
Further related work includes the literature on \emph{highest density
region} or \emph{level set} estimation [\citet{Polonik1995},
\citet{RigolletVert2009},
\citet{SamworthWand2010}].

\citet{HallSamworth2005} and \citet{BiauDevroye2010} proved an
analogous result for the bagged nearest neighbour classifier to the
\citet{Stone1977} result described in the previous paragraph. More
precisely, if the resample size $m = m_n$ used for the bagging diverges
to infinity, and $m/n \rightarrow0$ as $n \rightarrow\infty$, then
the risk of the bagged nearest neighbour classifier converges to the
Bayes risk. Note that this result does not depend on whether the
resamples are taken with or without replacement from the training data.
\citet{BCG2010} have recently proved a striking rate of convergence
result for the bagged nearest neighbour estimate; this is described in
greater detail in Section~\ref{SecBagged}.

\section{Main results}
\label{SecMain}

Let $(X,Y), (X_1,Y_1), (X_2,Y_2),\ldots$ be independent and identically
distributed pairs taking values in $\mathbb{R}^d \times\{1,2\}$. We
suppose that $\mathbb{P}(Y = 1) = \pi= 1- \mathbb{P}(Y=2)$ for some
$\pi\in(0,1)$ and that\vadjust{\goodbreak} $(X|Y=r) \sim P_r$ for $r= 1,2$, where~$P_r$
is a probability measure on $\mathbb{R}^d$. We write $\bar{P} = \pi P_1
+ (1-\pi)P_2$ for the marginal distribution of $X$ and let $\eta(x) =
\mathbb{P}(Y=1|X=x)$ denote the corresponding regression function.

A \emph{classifier} $C$ is a Borel measurable function from $\mathbb
{R}^d$ to $\{1,2\}$, with the interpretation that the point $x \in
\mathbb{R}^d$ is classified as belonging to class $C(x)$. The \emph
{misclassification rate}, or \emph{risk} of $C$ over a Borel measurable
set $\mathcal{R} \subseteq\mathbb{R}^d$ is defined to be
\[
R_{\mathcal{R}}(C) =
\mathbb{P}\bigl[\bigl\{C(X) \neq Y\bigr\} \cap \{X \in
\mathcal{R}\}\bigr].
\]
We also write $R(C)$ for this quantity when $\mathcal{R} = \mathbb
{R}^d$. The classifier which minimises the risk over $\mathcal{R}$ is
the Bayes classifier, given by
\[
C^{\mathrm{Bayes}}(x) = \cases{
1, & \quad $\mbox{if $\eta(x)
\geq1/2$},$
\vspace*{2pt}\cr
2, &\quad $ \mbox{otherwise.}$ }
\]
Its risk is
\[
R_{\mathcal{R}}\bigl(C^{\mathrm{Bayes}}\bigr) = \int_{\mathcal{R}}
\min\bigl\{\eta (x),1-\eta(x)\bigr\} \,d\bar{P}(x).
\]
For each $n \in\mathbb{N}$, let $\mathbf{w}_n = (w_{ni})_{i=1}^n$
denote a vector of weights, normalised so that $\sum_{i=1}^n w_{ni} =
1$. Fix $x \in\mathcal{R}$ and an arbitrary norm $\|\cdot\|$ on
$\mathbb{R}^d$, and let $(X_{(1)},Y_{(1)}),\ldots,(X_{(n)},Y_{(n)})$
denote a permutation of the training sample $(X_1,Y_1),\ldots
,(X_n,Y_n)$ such that $\|X_{(1)} - x\| \leq\cdots\leq\|X_{(n)} - x\|
$. We define the \emph{weighted nearest neighbour classifier} to be
\[
\hat{C}_n^{\mathrm{wnn}}(x) = \cases{ %
1, &\quad
$\mbox{if $\displaystyle\sum_{i=1}^n w_{ni}
\mathbh{1}_{\{
Y_{(i)}=1\}} \geq1/2$},$
\vspace*{2pt}\cr
2, & \quad$\mbox{otherwise.}$}
\]
We also write $\hat{C}_{n,\mathbf{w}_n}^{\mathrm{wnn}}$ where it is
necessary to emphasise the weight vector, for example, when comparing
different weighted nearest neighbour classifiers. Our initial goal is to
study the asymptotic behaviour of
\[
R_{\mathcal{R}}\bigl(\hat{C}_n^{\mathrm{wnn}}\bigr) = \mathbb{P}
\bigl[\bigl\{\hat {C}_n^{\mathrm{wnn}}(X) \neq Y\bigr\}
\mathbh{1}_{\{X \in\mathcal{R}\}} \bigr].
\]

It will be convenient to define some notation: for a smooth function
$g\dvtx \mathbb{R}^d \rightarrow\mathbb{R}$, we write $\dot{g}(x)$ for its
gradient vector at $x$, and $g_j(x)$ for its $j$th partial derivative
at $x$. Analogously, we write $\ddot{g}(x)$ for the Hessian matrix of
$g$ at $x$, and $g_{jk}(x)$ for its $(j,k)$th element. We let $B_\delta
(x) = \{y \in\mathbb{R}^d\dvtx \|y-x\| \leq\delta\}$ denote the closed
ball of radius $\delta$ centered at $x$ in the norm $\|\cdot\|$, and
let $a_d$ denote the $d$-dimensional Lebesgue measure of the unit ball
$B_1(x)$. Thus, $a_d = 2^d \Gamma(1+1/p)^d/\Gamma(1+d/p)$ when $\|\cdot
\|$ is the $\ell_p$-norm. We will make use of the following assumptions
for our theoretical results:
\begin{longlist}[(A.1)]
\item[(A.1)] The set $\mathcal{R} \subseteq\mathbb{R}^d$ is a compact
$d$-dimensional manifold with boundary~$\partial\mathcal{R}$.
\item[(A.2)] The set $\mathcal{S} = \{x \in\mathcal{R}\dvtx \eta(x) = 1/2\}
$ is nonempty. There exists an open subset $U_0$ of $\mathbb{R}^d$
that contains $\mathcal{S}$ and such that the following properties
hold: first,
$\eta$ is continuous on $U \setminus U_0$,
where $U$ is an open set containing $\mathcal{R}$;
second, the restrictions of $P_1$ and $P_2$ to $U_0$ are absolutely
continuous with respect to Lebesgue measure, with twice continuously
differentiable Radon--Nikodym derivatives $f_1$ and $f_2$, respectively.
%
\item[(A.3)] There exists $\rho> 0$ such that $\int_{\mathbb{R}^d} \|
x\|^\rho \,d\bar{P}(x) < \infty$. Moreover, for sufficiently small
$\delta> 0$, the ratio $\bar{P}(B_\delta(x))/(a_d \delta^d)$ is
bounded away from zero, uniformly for $x \in\mathcal{R}$.
\item[(A.4)] For all $x \in\mathcal{S}$, we have $\dot{\eta}(x) \neq
0$, and for all $x \in\mathcal{S} \cap\partial\mathcal{R}$, we have
$\dot{\partial\eta}(x) \neq0$, where $\partial\eta$ denotes the
restriction of $\eta$ to $\partial\mathcal{R}$.
\end{longlist}
The introduction of the compact set $\mathcal{R}$ finesses the problem
of performing classification in the tails of the feature vector
distributions. See, for example, \citet{HallKang2005}, Section~3, for
further discussion of this point and related results, as well as \citet
{ChandaRuymgaart1989}. \citet{MammenTsybakov1999} and \citet
{AudibertTsybakov2007} impose similar compactness assumptions for their
results. The set $\mathcal{R}$ may be arbitrarily large, though the
larger it is, the stronger are the requirements in ({A.2}).
Although as stated, the assumptions on $\mathcal{R}$ are quite general,
little is lost by thinking of $\mathcal{R}$ as a large closed Euclidean
ball. Its role in the asymptotic expansion of Theorem~\ref
{ThmOptWeights} below is that it is involved in the definition of the
set $\mathcal{S}$, which represents the decision boundary of the Bayes
classifier. We will see that the behaviour of $f_1$ and $f_2$ on the set
$\mathcal{S}$ is crucial for determining the asymptotic behaviour of
weighted nearest neighbour classifiers.

The second part of ({A.3}) asks that the ratio of the $\bar
{P}$-measure of small balls to the corresponding $d$-dimensional
Lebesgue measure is bounded away from zero. This requirement is
satisfied, for instance, if $P_1$ and $P_2$ are absolutely continuous
with respect to Lebesgue measure, with Radon--Nikodym derivatives that
are bounded away from zero on the open set $U$.

The assumption in ({A.4}) that $\dot{\eta}(x) \neq0$ for $x \in
\mathcal{S}$ asks that $f_1$ and $f_2$, weighted by the respective
prior probabilities of each class, should cut at a nonzero angle along
$\mathcal{S}$. In the language of differential topology, this means
that $1/2$ is a \emph{regular value} of the function $\eta$, and the
second part of ({A.4}) asks for $1/2$ to be a regular value of the
restriction of $\eta$ to $\partial\mathcal{R}$. Together, these two
requirements ensure that $\mathcal{S}$ is a $(d-1)$-dimensional
submanifold with boundary of $\mathbb{R}^d$, and the boundary of
$\mathcal{S}$ is $\{x \in\partial\mathcal{R}\dvtx \eta(x) = 1/2\}$
[\citet{GuilleminPollack1974}, page 60].

The requirement in~({A.4}) that $\dot{\eta}(x) \neq0$ for $x
\in\mathcal{S}$ is related to the well-known \emph{margin condition}
of, for example, \citet{MammenTsybakov1999} and \citet{Tsybakov2004};
when it holds (and in the presence of the other conditions), there
exist $c,C > 0$ such that
%
\begin{equation}
\label{EqMargin} c \varepsilon\leq\mathbb{P}\bigl(\bigl|\eta(X) - 1/2\bigr| \leq
\varepsilon\cap X \in \mathcal{R}\bigr) \leq C \varepsilon
\end{equation}
for sufficiently small $\varepsilon> 0$; see \citet{Tsybakov2004}, Proposition~1. A proof of this fact, which uses Weyl's
tube formula [\citet{Gray2004}], is given after the completion of the
proof of Theorem~\ref{ThmMain} in the supplementary material
[\citet{Samworth2012}]. In this sense, we work in the setting of a margin
condition with the power parameter equal to 1.

We now introduce some notation needed for Theorem~\ref{ThmMain} below.
For $\beta> 0$, let $W_{n,\beta}$ denote the set of all sequences of
nonnegative deterministic weight vectors $\mathbf{w}_n =
(w_{ni})_{i=1}^n$ satisfying:
\begin{itemize}
\item$\sum_{i=1}^n w_{ni}^2 \leq n^{-\beta}$;
\item$n^{-4/d} (\sum_{i=1}^n \alpha_i w_{ni} )^2 \leq n^{-\beta
}$, where $\alpha_i = i^{1+2/d} - (i-1)^{1+2/d}$; note that this latter
expression appears in~(\ref{EqOptWeights});
\item$n^{2/d}\sum_{i=k_2+1}^n w_{ni}/\sum_{i=1}^n \alpha_i w_{ni} \leq
1/\log n$, where $k_2 = \lfloor n^{1-\beta} \rfloor$;
\item$\sum_{i=k_2+1}^n w_{ni}^2/\sum_{i=1}^n w_{ni}^2 \leq1/\log n$;
\item$\sum_{i=1}^n w_{ni}^3/(\sum_{i=1}^n w_{ni}^2)^{3/2} \leq1/\log n$.
\end{itemize}
Observe that $W_{n,\beta_1} \supset W_{n,\beta_2}$ for $\beta_1 < \beta_2$.
The first and last conditions ensure that the weights are not too
concentrated on a small number of points; the second amounts to a mild
moment condition on the probability distribution on $\{1,\ldots,n\}$
defined by the weights. The next two conditions ensure that not too
much weight (or squared weight in the case of the latter condition) is
assigned to observations that are too far from the point being
classified. Although there are many requirements on the weight vectors,
they are rather mild conditions when $\beta$ is small, as can be seen
by considering the limiting case $\beta= 0$. For instance, for the
unweighted $k$-nearest neighbour classifier with weights $\mathbf{w}_n =
(w_{ni})_{i=1}^n$ given by $w_{ni}=k^{-1} \mathbh{1}_{\{1 \leq i \leq
k\}}$, we have that $\mathbf{w}_n \in W_{n,\beta}$ for small $\beta>
0$ provided that $\max(n^{\beta},\log^2 n) \leq k \leq\min(n^{(1-\beta
d/4)},n^{1-\beta})$. Thus for the vector of $k$-nearest neighbour
weights to belong to $W_{n,\beta}$ for all large $n$, it is necessary
that the usual conditions $k \rightarrow\infty$ and $k/n \rightarrow
0$ for consistency are satisfied, and these conditions are almost
sufficient when $\beta> 0$ is small. The situation is similar for the
bagged nearest neighbour classifier---see Section~\ref{SecBagged} below.

The fact that the weights are assumed to be deterministic means that
they depend only on the ordering of the distances, not the raw
distances themselves (as would be the case for a classifier based on
kernel density estimates of the population densities). Such
kernel-based classifiers are not necessarily straightforward to
implement, however: \citet{HallKang2005} showed that even in the simple
situation where $d=1$ and $\pi f_1$ and $(1-\pi)f_2$ cross at a single
point $x_0$, the optimal order of the bandwidth for the kernel depends
on the sign of $\ddot{f}_1(x_0)\ddot{f}_2(x_0)$.

Continuing with our notational definitions, let $\bar{f} = \pi f_1 +
(1-\pi)f_2$. Define
%
\begin{equation}
\label{Eqa} a(x) = \frac{\sum_{j=1}^d c_{j,d}\{\eta_j(x)\bar{f}_j(x) +
({1}/{2})\eta_{jj}(x)\bar{f}(x)\}}{a_d^{1+2/d}\bar{f}(x)^{1+2/d}},
\end{equation}
where $c_{j,d} = \int_{v:\|v\| \leq1} v_j^2 \,dv$. Finally, let
%
\begin{eqnarray}
\label{EqC1C2} B_1& = &\int_{\mathcal{S}}
\frac{\bar{f}(x_0)}{4\|\dot{\eta}(x_0)\|} \,d\mathrm{Vol}^{d-1}(x_0)
\quad\mbox{and}
\nonumber
\\[-8pt]
\\[-8pt]
\nonumber
B_2 &=& \int_{\mathcal{S}} \frac{\bar
{f}(x_0)}{\|\dot{\eta}(x_0)\|}a(x_0)^2
\,d\mathrm{Vol}^{d-1}(x_0),
\end{eqnarray}
where $\mathrm{Vol}^{d-1}$ denotes the natural $(d-1)$-dimensional
volume measure that $\mathcal{S}$ inherits as a subset of $\mathbb
{R}^d$. Note that $B_1 > 0$, and $B_2 \geq0$, with equality if and
only if~$a$ is identically zero on $\mathcal{S}$. Although the
definitions of $B_1$ and $B_2$ are complicated, we will see after the
statement of Theorem~\ref{ThmMain} below that they are comprised of
terms that have natural interpretations.
%
\begin{thmm}
\label{ThmMain}
Assume \textup{({A.1})}, \textup{({A.2})}, \textup{({A.3})} and \textup{({A.4})}.
Then for each $\beta\in(0,1/2)$,
\[
R_{\mathcal{R}}\bigl(\hat{C}_n^{\mathrm{wnn}}\bigr) -
R_{\mathcal{R}}\bigl(C^{\mathrm
{Bayes}}\bigr) = \gamma_n(
\mathbf{w}_n)\bigl\{1 + o(1)\bigr\}
\]
as $n \rightarrow\infty$, uniformly for $\mathbf{w}_n \in W_{n,\beta
}$, where
\[
\gamma_n(\mathbf{w}_n) = B_1 \sum
_{i=1}^n w_{ni}^2 +
B_2 \Biggl(\sum_{i=1}^n
\frac{\alpha_i w_{ni}}{n^{2/d}} \Biggr)^2.
\]
\end{thmm}
Theorem~\ref{ThmMain} tells us that, asymptotically, the dominant
contribution to the regret over $\mathcal{R}$ of the weighted nearest
neighbour classifier can be decomposed as a sum of two terms. The two
terms, constant multiples of $\sum_{i=1}^n w_{ni}^2$ and $ (\sum_{i=1}^n \frac{\alpha_i w_{ni}}{n^{2/d}} )^2$, respectively,
represent variance and squared bias contributions to the regret. It is
interesting to observe that, although the 0--1 classification loss
function is quite different from the squared error loss often used in
regression problems, we nevertheless obtain such an asymptotic decomposition.

The constant multiples of the dominant variance and squared bias terms
depend only on the behaviour of $f_1$ and $f_2$ (and their first and
second derivatives) on $\mathcal{S}$, as seen from~(\ref{EqC1C2}).
Moreover, we can see from the expression for $B_1$ in~(\ref{EqC1C2})
that the contribution to the dominant variance term in the regret will
tend to be large in the following three situations: first, when $\bar
{f}(\cdot)$ is large on $\mathcal{S}$; second, when the $\mathrm
{Vol}^{d-1}$ measure of $\mathcal{S}$ is large; and third, when $\|\dot
{\eta}(\cdot)\|$ is small on $\mathcal{S}$. In the first two of these
situations, the probability is relatively high that a point to be
classified will be close to the Bayes decision boundary $\mathcal{S}$,
where classification is difficult. In the latter case, the regression
function $\eta$ moves away from $1/2$ only slowly as we move away from~$\mathcal{S}$,
meaning that there is a relatively large region of
points near $\mathcal{S}$ where classification is difficult. From the
expression for $B_2$ in~(\ref{EqC1C2}), we see that the dominant
squared bias term is also large in these situations, and also when
$a(\cdot)^2$ is large on $\mathcal{S}$. From the proof of Theorem~\ref
{ThmMain}, it is apparent that $a(x)\sum_{i=1}^n \frac{\alpha_i
w_{ni}}{n^{2/d}}$ is the dominant bias term for $S_n(x) = \sum_{i=1}^n
w_{ni} \mathbh{1}_{\{Y_{(i)} = 1\}}$ as an estimator of $\eta(x)$.
Indeed, by a Taylor expansion,
\begin{eqnarray*}
&&\mathbb{E}\bigl\{S_n(x)\bigr\} - \eta(x) \\
&&\qquad= \sum
_{i=1}^n w_{ni}\mathbb{E}\eta
(X_{(i)}) - \eta(x)
\\
&&\qquad\approx\sum_{i=1}^n w_{ni}
\mathbb{E}\bigl\{(X_{(i)}-x)^T\dot{\eta}(x)\bigr\} +
\frac{1}{2}\sum_{i=1}^n
w_{ni} \mathbb{E}\bigl\{(X_{(i)}-x)^T\ddot{\eta
}(x) (X_{(i)}-x)\bigr\}.
\end{eqnarray*}
The two summands in the definition of $a(x)$ represent asymptotic
approximations to the respective summands in this approximation.

Consider now the problem of optimising the choice of weight vectors. Let
%
\begin{equation}
\label{Eqkstar} k^* = \biggl\lfloor \biggl\{\frac{d(d+4)}{2(d+2)} \biggr
\}^{d/(d+4)} \biggl(\frac{B_1}{B_2} \biggr)^{d/(d+4)}
n^{4/(d+4)} \biggr\rfloor,
\end{equation}
and then define the weights $\mathbf{w}_n^* = (w_{ni}^*)_{i=1}^n$ as
in~(\ref{EqOptWeights}). The first part of Theorem~\ref
{ThmOptWeights} below can be regarded as saying that the weights
$\mathbf{w}_n^*$ are asymptotically optimal.
%
\begin{thmm}
\label{ThmOptWeights}
Assume \textup{{(A.1)}--({A.4})}, and assume also that $B_2 > 0$.
For any $\beta> 0$ and any sequence $\mathbf{w}_n = (w_{ni})_{i=1}^n
\in W_{n,\beta}$, we have
%
\begin{equation}
\label{EqRegRatio} \liminf_{n \rightarrow\infty} \frac{R_{\mathcal{R}}(\hat{C}_{n,\mathbf
{w}_n}^{\mathrm{wnn}}) - R_{\mathcal{R}}(C^{\mathrm
{Bayes}})}{R_{\mathcal{R}}(\hat{C}_{n,\mathbf{w}_n^*}^{\mathrm{wnn}}) -
R_{\mathcal{R}}(C^{\mathrm{Bayes}})} \geq1.
\end{equation}
Moreover, the ratio in~(\ref{EqRegRatio}) above converges to 1 if and
only if we have both $\sum_{i=1}^n w_{ni}^2/\sum_{i=1}^n (w_{ni}^*)^2
\rightarrow1$ and $\sum_{i=1}^n \alpha_i w_{ni}/\sum_{i=1}^n \alpha_i
w_{ni}^* \rightarrow1$. Equivalently, this occurs if and only if both
%
\begin{eqnarray}
\label{EqEquivCond} n^{4/(d+4)}\sum_{i=1}^n
\bigl\{w_{ni}^2 - \bigl(w_{ni}^*
\bigr)^2\bigr\} &\rightarrow&0\quad \mbox{and}
\nonumber
\\[-8pt]
\\[-8pt]
\nonumber
 n^{-{8}/{(d(d+4))}} \sum
_{i=1}^n \alpha_i
\bigl(w_{ni} - w_{ni}^*\bigr) &\rightarrow&0.
\end{eqnarray}
Finally,
%
\begin{eqnarray}
\label{EqAsympRegret}&& n^{4/(d+4)} \bigl\{R_{\mathcal{R}}\bigl(
\hat{C}_{n,\mathbf{w}_n^*}^{\mathrm
{wnn}}\bigr) - R_{\mathcal{R}}
\bigl(C^{\mathrm{Bayes}}\bigr)\bigr\}
\nonumber
\\[-8pt]
\\[-8pt]
\nonumber
&&\qquad \rightarrow\frac
{(d+2)^{{(2d+4)}/{(d+4)}}}{2^{{4}/{(d+4)}}} \biggl(
\frac{d+4}{d} \biggr)^{{d}/{(d+4)}}B_1^{{4}/{(d+4)}}B_2^{{d}/{(d+4)}}.
\end{eqnarray}
\end{thmm}
Now write $\hat{C}_{n,k}^{\mathrm{nn}}$ for the traditional, unweighted
$k$-nearest neighbour classifier (or equivalently, the weighted nearest
neighbour classifier with $w_{ni} = 1/k$ for $i=1,\ldots,k$ and $w_{ni}
= 0$ otherwise). Another consequence of Theorem~\ref{ThmMain} is that,
provided {(A.1)}--({A.4}) hold and $B_2 > 0$, the
quantity $k^*$ defined in~(\ref{Eqkstar}) is larger by a factor of
$ \{\frac{2(d+4)}{d+2} \}^{d/(d+4)}$ (up to an unimportant
rounding error) than the asymptotically optimal choice of $k^{\mathrm
{opt}}$ for $\hat{C}_{n,k}^{\mathrm{nn}}$; see also \citet{HPS2008}. We
can therefore compare the performance of $\hat{C}_{n,k^{\mathrm
{opt}}}^{\mathrm{nn}}$ with that of $\hat{C}_{n,\mathbf{w}_n^*}^{\mathrm{wnn}}$.
%
\begin{cor}
\label{CorRegretRatio}
Assume \textup{{(A.1)}--({A.4})} and assume also that $B_2 > 0$. Then
%
\begin{equation}
\label{EqRegretRatio} \frac{R_{\mathcal{R}}(\hat{C}_{n,\mathbf{w}_n^*}^{\mathrm{wnn}}) -
R_{\mathcal{R}}(C^{\mathrm{Bayes}})}{R_{\mathcal{R}}(\hat
{C}_{n,k^{\mathrm{opt}}}^{\mathrm{nn}}) - R_{\mathcal{R}}(C^{\mathrm
{Bayes}})} \rightarrow\frac{1}{4^{d/(d+4)}} \biggl(
\frac
{2d+4}{d+4} \biggr)^{(2d+4)/(d+4)}
\end{equation}
as $n \rightarrow\infty$.
\end{cor}
Since the limit in~(\ref{EqRegretRatio}) does not depend on the
underlying populations, we can plot it as a function of $d$; cf.
Figure~\ref{FigAsympImp}. In fact, Corollary~\ref{CorRegretRatio}
suggests a natural correspondence between any unweighted $k$-nearest
neighbour classifier $\hat{C}_{n,k}^{\mathrm{nn}}$ and the weighted
nearest neighbour classifier which we denote by $\hat{C}_{n,\mathbf
{w}_n^{\mu(k)}}^{\mathrm{wnn}}$ whose weights are of the optimal
form~(\ref{EqOptWeights}), but with $k^*$ replaced with
%
\begin{equation}
\label{Eqmuk} \mu(k) = \biggl\lfloor \biggl\{\frac{2(d+4)}{d+2} \biggr
\}^{d/(d+4)}k \biggr\rfloor.
\end{equation}
Under the conditions of Corollary~\ref{CorRegretRatio}, we can compare
$\hat{C}_{n,k}^{\mathrm{nn}}$ and $\hat{C}_{n,\mathbf{w}_n^{\mu
(k)}}^{\mathrm{wnn}}$, concluding that for each $\beta\in(0,1/2)$,
%
\begin{equation}
\label{EqRegretRatio2} \frac{R_{\mathcal{R}}(\hat{C}_{n,\mathbf{w}_n^{\mu(k)}}^{\mathrm{wnn}})
- R_{\mathcal{R}}(C^{\mathrm{Bayes}})}{R_{\mathcal{R}}(\hat
{C}_{n,k}^{\mathrm{nn}}) - R_{\mathcal{R}}(C^{\mathrm{Bayes}})} \rightarrow\frac{1}{4^{d/(d+4)}} \biggl(
\frac{2d+4}{d+4} \biggr)^{(2d+4)/(d+4)}
\end{equation}
as $n \rightarrow\infty$, uniformly for $n^{\beta} \leq k \leq
n^{1-\beta}$. The fact that the convergence in~(\ref{EqRegretRatio2})
is uniform for $k$ in this range means that the ratio on the left-hand
side of~(\ref{EqRegretRatio2}) has the same limit if we replace $k$ by
an estimator $\hat{k}$ constructed from the training data
$(X_1,Y_1),\ldots,(X_n,Y_n)$, provided that $\hat{k}$ lies in this
range with probability tending to 1.

In a complementary approach to that taken in most of this paper, \citet
{AudibertTsybakov2007} study the minimax properties of plug-in
classifiers. They show in particular that a certain classifier obtained
by modifying a local polynomial estimator of the regression function
$\eta$ attains the minimax rate over a set of distributions $P$ of
random vectors $(X,Y)$ on $\mathbb{R}^d \times\{1,2\}$ for which the
regression function belongs to a H\"older class, $P$\vadjust{\goodbreak} satisfies a margin
condition and the marginal distribution of $X$ satisfies a so-called
\emph{strong density assumption}. This rate is $O(n^{-4/(d+4)})$ when
the H\"older smoothness parameter is 2, and the margin power parameter
is 1. By adapting their arguments, we are able to show in the
supplementary material [\citet{Samworth2012}] that several weighted
nearest-neighbour classifiers (including the unweighted, optimally
weighted and bagged versions of Section~\ref{SecBagged}) can also
attain this minimax rate. Such results give reassurance about
worst-case behaviour; however, they do not lead naturally to an optimal
weighting scheme or a quantification of the relative performance of two
weighted nearest neighbour classifiers attaining the same rate. These
are the main goals of this work.

Finally in this section, we note that the theory presented above can be
extended in a natural way to multicategory classification problems,
where the class labels take values in the set $\{1,\ldots,K\}$. Writing
$\eta_r(x) = \mathbb{P}(Y=r|X=x)$, let
\[
\mathcal{S}_{r_1,r_2} = \Bigl\{x \in\mathcal{R}\dvtx \mathop{\argmax}_{r \in\{
1,\ldots,K\}}
\eta_r(x) = \{r_1,r_2\} \Bigr\}
\]
for distinct indices $r_1,r_2 \in\{1,\ldots,K\}$. In addition to
{(A.1)} and the obvious analogues of the conditions
{(A.2)}, {(A.3)} and~({A.4}), we require:
\begin{longlist}[(A.5)]
\item[(A.5)] For each $(r_1,r_2) \neq(r_3,r_4)$, the submanifolds
$\mathcal{S}_{r_1,r_2}$ and $\mathcal{S}_{r_3,r_4}$ of $\mathbb{R}^d$
are transversal.
\end{longlist}
Condition ({A.5}) ensures that $\mathcal{S}_{r_1,r_2} \cap
\mathcal{S}_{r_3,r_4} \cap(\mathcal{R} \setminus\partial\mathcal
{R})$ is either empty or a $(d-2)$-dimensional submanifold of $\mathbb
{R}^d$ [\citet{GuilleminPollack1974}, page 30]. Under these conditions,
the conclusion of Theorem~\ref{ThmMain} holds, provided that the
constants~$B_1$ and $B_2$ are replaced with $\tilde{B}_1 = \sum_{r_1
\neq r_2} B_{1,r_1,r_2}$ and $\tilde{B}_2 = \sum_{r_1 \neq r_2}
B_{2,r_1,r_2}$, respectively, where each term $B_{1,r_1,r_2}$ and
$B_{2,r_1,r_2}$ is an integral over $\mathcal{S}_{r_1,r_2}$. Apart from
the obvious notational changes involved in converting $B_1$ and $B_2$
to $B_{1,r_1,r_2}$ and $B_{2,r_1,r_2}$, the only other change required
is to replace the constant factor $1/4$ in the definition of $B_1$ with
$\eta_{r_1,r_2}(x_0)\{1 - \eta_{r_1,r_2}(x_0)\}$ where $\eta_{r_1,r_2}(x_0)$
denotes the common value that $\eta_{r_1}$ and $\eta_{r_2}$ take at $x_0 \in\mathcal{S}_{r_1,r_2}$. This change accounts
for the fact that $\eta_{r_1,r_2}(x_0)$ is not necessarily equal to
$1/2$ on $\mathcal{S}_{r_1,r_2}$.

It follows (provided also that $\tilde{B}_2 > 0$) that the
asymptotically optimal weights are still of the form~(\ref
{EqOptWeights}), but with the ratio $B_1/B_2$ in the expression for
$k^*$ in~(\ref{Eqkstar}) replaced with $\tilde{B}_1/\tilde{B}_2$.
Moreover, the conclusion of Corollary~\ref{CorRegretRatio} and the
subsequent discussion also remain true.

\section{The bagged nearest neighbour classifier}
\label{SecBagged}

Traditionally, the bagged nearest neighbour classifier is obtained by
applying the 1-nearest neighbour classifier to many resamples from the
training data. The final classification is made by a majority vote on
the classifications obtained from the resamples. In the most common\vadjust{\goodbreak}
version of bagging where the resamples are drawn with replacement, and
the resample size is the same as the original sample size, bagging the
nearest neighbour classifier gives no improvement over the 1-nearest
neighbour classifier [\citet{HallSamworth2005}]. This is because the
nearest neighbour occurs in more than half (in fact, roughly a
proportion $1-1/e$) of the resamples.

Nevertheless, if a smaller resample size is used, then substantial
improvements over the nearest neighbour classifier are possible, as has
been verified empirically by \citet{MartinezMunozSuarez2010}. In fact,
if the resample size is $m$, then the ``infinite simulation'' versions
of the bagged nearest neighbour classifier in the with- and
without-replacement resampling cases are weighted nearest neighbour
classifiers with respective weights
%
\begin{equation}
\label{Eqbagged1} w_{ni}^{\mathrm{b},\mathrm{with}} = \biggl(1 -
\frac{i-1}{n} \biggr)^m - \biggl(1 - \frac{i}{n}
\biggr)^m, \qquad i=1,\ldots,n
\end{equation}
and
%
\begin{equation}
\label{Eqbagged2} w_{ni}^{\mathrm{b},\mathrm{w/o}} = \cases{ \pmatrix{n-i
\cr
m-1}\Big/\pmatrix {n
\cr
m}, &\quad $\mbox{for $i=1,\ldots,n-m+1$},$
\vspace*{2pt}\cr
0, & \quad $\mbox{for
$i=n-m+2,\ldots,n$.}$}
\end{equation}
Of course, the observations above render the resampling redundant, and
we regard the weighted nearest neighbour classifiers with the weights
above as defining the two versions of the bagged nearest neighbour
classifier. It is convenient to let $q = m/n$ denote the resampling
fraction. Intuitively, for large $n$, both versions of the bagged
nearest neighbour classifier behave like the weighted nearest neighbour
classifier with weights $(w_{ni}^{\mathrm{Geo}})_{i=1}^n$ which place a
$\mathrm{Geometric}(q)$ distribution (conditioned on being in the set
$\{1,\ldots,n\}$) on the weights
%
\begin{equation}
\label{Eqbagged3} w_{ni}^{\mathrm{Geo}} = \frac{q(1-q)^{i-1}}{1 - (1-q)^n},\qquad i=1,
\ldots,n.
\end{equation}
The reason for this is that, in order for the $i$th nearest neighbour of
the training data to be the nearest neighbour of the resample, the
nearest $i-1$ neighbours must not appear in the resample, while the
$i$th nearest neighbour must appear, and these events are almost
independent when $n$ is large; see \citet{HallSamworth2005}. Naturally,
the parameter $q$ plays a crucial role in the performance of the bagged
nearest neighbour classifier, and for small $\beta> 0$, the three
vectors of weights given in~(\ref{Eqbagged1}), (\ref{Eqbagged2})
and~(\ref{Eqbagged3}) belong to $W_{n,\beta}$ for all large $n$ if
$\max(\frac{1}{2}n^{-(1-\beta d/4)},n^{-(1-2\beta)}) \leq q \leq
3n^{-\beta}$. In the following corollary of Theorem~\ref{ThmMain}, we
write $\hat{C}_{n,q}^{\mathrm{bnn}}$ to denote either of the bagged
nearest neighbour classifiers with weights~(\ref{Eqbagged1}),~(\ref
{Eqbagged2}) or their approximation with weights~(\ref{Eqbagged3}).
%
\begin{cor}
\label{CorBaggednnExp}
Assume \textup{({A.1})--({A.4})}. For every $\beta\in(0,1/2)$,
\[
R_{\mathcal{R}}\bigl(\hat{C}_{n,q}^{\mathrm{bnn}}\bigr) -
R_{\mathcal
{R}}\bigl(C^{\mathrm{Bayes}}\bigr) = \tilde{\gamma}_n(q)
\bigl\{1+o(1)\bigr\},\vadjust{\goodbreak}
\]
uniformly for $n^{-(1-\beta)} \leq q \leq n^{-\beta}$, where
\[
\tilde{\gamma}_n(q) = \frac{B_1}{2}q + \frac{B_2 \Gamma (2 +
{2}/{d} )^2}{n^{4/d}q^{4/d}}.
\]
\end{cor}
This result is somewhat related to Corollary~10 of \citet{BCG2010}. In
that paper, the authors study the bagged nearest neighbour estimate $\hat
{\eta}_n$ of the regression function $\eta$. They prove in particular
that under regularity conditions (including a Lipschitz assumption on
$\eta$) and for a suitable choice of resample size,
\[
\mathbb{E} \bigl[\bigl\{\hat{\eta}_n(X) - \eta(X)\bigr
\}^2 \bigr] = O\bigl(n^{-2/(d+2)}\bigr)
\]
for $d \geq3$. It is known [e.g., Ibragimov and Khasminski{\u\i}
(\citeyear{IbragimovKhasminskii1980,IbragimovKhasminskii1981,IbragimovKhasminskii1982})]
that this is the minimax optimal rate for their problem.

\begin{figure}

\includegraphics{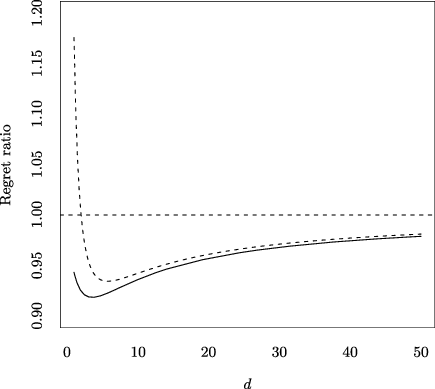}

\caption{Asymptotic ratio of the regret of the bagged
nearest neighbour classifier (dashed) to that of the $k$-nearest
neighbour classifier,
as a function of the dimension of the feature vectors. The asymptotic
regret ratio
for the optimally weighted nearest neighbour classifier compared with
the $k$-nearest
neighbour classifier is shown as a solid line for comparison.}\label
{FigRegretRatio2}
\end{figure}

Corollary~\ref{CorBaggednnExp} may also be applied to deduce that the
asymptotically optimal choice of $q$ in all three cases is
\[
q^{\mathrm{opt}} = \frac{8^{d/(d+4)}\Gamma (2+{2}/{d}
)^{2d/(d+4)}}{d^{d/(d+4)}} \biggl(\frac{B_2}{B_1}
\biggr)^{d/(d+4)} n^{-4/(d+4)}.
\]
Thus, in an analogous fashion to Section~\ref{SecMain}, we can
consider the performance of $\hat{C}_{n,q^{\mathrm{opt}}}^{\mathrm
{bnn}}$ relative to that of $\hat{C}_{n,k^{\mathrm{opt}}}^{\mathrm{nn}}$.
%
\begin{cor}
\label{CorBagged}
Assume \textup{{(A.1)}--({A.4})} and assume also that $B_2 > 0$. Then
%
\begin{equation}
\label{EqRegretRatio3} \frac{R_{\mathcal{R}}(\hat{C}_{n,q^{\mathrm{opt}}}^{\mathrm{bnn}}) -
R_{\mathcal{R}}(C^{\mathrm{Bayes}})}{R_{\mathcal{R}}(\hat
{C}_{n,k^{\mathrm{opt}}}^{\mathrm{nn}}) - R_{\mathcal{R}}(C^{\mathrm
{Bayes}})} \rightarrow\frac{\Gamma (2+{2}/{d}
)^{2d/(d+4)}}{2^{4/(d+4)}}
\end{equation}
as $n \rightarrow\infty$.
\end{cor}
The limiting ratio in~(\ref{EqRegretRatio3}) is plotted as a function
of $d$ in Figure~\ref{FigRegretRatio2}. The ratio is about 1.18 when
$d=1$, showing that the bagged nearest neighbour classifier has
asymptotically worse performance than the $k$-nearest neighbour
classifier in this case. The ratio is equal to 1 when $d=2$, and is
less than 1 for $d \geq3$. The facts that the asymptotically optimal
weights decay as illustrated in Figure~\ref{FigDecWeights} and that
the bagged nearest neighbour weights decay approximately geometrically
explain why the bagged nearest neighbour classifier has almost optimal
performance among nonnegatively weighted nearest neighbour classifiers
when $d$ is large.

Similar to the discussion following Corollary~\ref{CorRegretRatio},
based on the expressions for $k^{\mathrm{opt}}$ and $q^{\mathrm{opt}}$,
there is a natural correspondence between the unweighted $k$-nearest
neighbour classifier $\hat{C}_{n,\hat{k}}^{\mathrm{nn}}$ with data
driven $\hat{k}$, and the bagged nearest neighbour classifier $\hat
{C}_{n,\hat{q}}^{\mathrm{bnn}}$, where
%
\begin{equation}
\label{Eqhatq} \hat{q} = 2^{d/(d+4)}\Gamma \biggl(2+\frac{2}{d}
\biggr)^{2d/(d+4)} \frac
{1}{\hat{k}}.
\end{equation}
The same limit~(\ref{EqRegretRatio3}) holds for the regret ratio of
these classifiers, again provided there exists $\beta\in(0,1/2)$ such
that $\mathbb{P}(n^{\beta} \leq\hat{k} \leq n^{1-\beta}) \rightarrow1$.

%

\section{Faster rates of convergence}
\label{SecNegWeights}

If we allow negative weights, it is possible to choose weights
satisfying $\sum_{i=1}^n \alpha_i w_{ni} = 0$. This means that we can
eradicate the dominant squared bias term in the asymptotic expansion of
Theorem~\ref{ThmMain}. It follows that, subject to additional
smoothness conditions, we can achieve faster rates of convergence with
weighted nearest neighbour classifiers, as we now describe. The
appropriate variant of condition~({A.2}), which we denote
by~{(A.2)(r)}, is as follows:
\begin{longlist}[(A.2)(r)]
\item[(A.2)(r)] The set $\mathcal{S} = \{x \in\mathcal{R}\dvtx \eta(x) =
1/2\}$ is nonempty. There exists an open subset $U_0$ of $\mathbb
{R}^d$ that contains $\mathcal{S}$ and such that the following
properties hold: first,
$\eta$ is continuous on $U \setminus U_0$,
where $U$ is an open set containing $\mathcal
{R}$; second, the restrictions of $P_1$ and $P_2$ to $U_0$ are
absolutely continuous with respect to Lebesgue measure, with $2r$-times
continuously differentiable Radon--Nikodym derivatives $f_1$ and $f_2$,
respectively.
\end{longlist}
Thus condition {(A.2)(1)} is identical to {(A.2)}. Note
that we are still in the setting of a margin condition with power
parameter equal to 1. Let $S$ denote the set of multi-indices $s =
(s_1,\ldots,s_d)$, so $s$ is a $d$-tuple of nonnegative integers. For
$s \in S$, we write $|s| = s_1 + \cdots+ s_d$, and for $v = (v_1,\ldots
,v_d)^T \in\mathbb{R}^d$, we write $v^s = v_1^{s_1}v_2^{s_2}\cdots
v_d^{s_d}$. Now, for $s \in S$, let $c_{s,d} = \int_{\|v\| \leq1} v^s
\,dv$. It is convenient here to use multi-index notation for derivatives,
so we write $g_s(x) = \frac{\partial^{|s|}}{\partial x_1^{s_1} \cdots
\partial x_d^{s_d}}g(x)$. Now let
\[
\bar{S}_r = \bigl\{\bigl(s^1,s^2\bigr) \in
S \times S\dvtx \bigl|s^1\bigr| + \bigl|s^2\bigr| = 2r, \bigl|s^1\bigr|
\geq 1, s_j^1 + s_j^2 \in2
\mathbb{Z}\ \forall j = 1,\ldots,d\bigr\},
\]
and let
\[
a^{(r)}(x) = \frac{1}{a_d^{1+2r/d}\bar{f}(x)^{1+2r/d}}\sum_{(s^1,s^2)
\in\bar{S}_r}
\frac{c_{s^1+s^2,d} \eta_{s^1}(x)\bar{f}_{s^2}(x)}{|s^1|!|s^2|!};
\]
thus $a^{(1)}(x) = a(x)$. Further, let
\[
B_2^{(r)} = \int_{\mathcal{S}}
\frac{\bar{f}(x_0)}{\|\dot{\eta}(x_0)\|
}a^{(r)}(x_0)^2 \,d
\mathrm{Vol}^{d-1}(x_0).
\]
For $\ell\in\mathbb{N}$, define $\alpha_i^{(\ell)} = i^{1+2\ell/d} -
(i-1)^{1+2\ell/d}$. We consider restrictions on the set of weight
vectors analogous to those imposed on $r$th order kernels in kernel
density estimation. Specifically, we let $W_{n,\beta,r}^\dagger$ denote
the set of deterministic weight vectors $\mathbf{w}_n =
(w_{ni})_{i=1}^n$ satisfying:
\begin{itemize}
\item$\sum_{i=1}^n w_{ni} = 1$, $n^{2r/d}\sum_{i=1}^n \alpha_i^{(\ell
)}w_{ni}/n^{2\ell/d}\sum_{i=1}^n \alpha_i^{(r)}w_{ni} \leq1/\log n$
for $\ell=\break 1,\ldots,r-1$;
\item$\sum_{i=1}^n w_{ni}^2 \leq n^{-\beta}$;
\item$n^{-4r/d} (\sum_{i=1}^n \alpha_i^{(r)} w_{ni} )^2 \leq
n^{-\beta}$;
%
\item there exists $k_2 \leq\lfloor n^{1-\beta} \rfloor$ such that
$n^{2r/d}\sum_{i=k_2+1}^n |w_{ni}|/\sum_{i=1}^n \alpha_i^{(r)} w_{ni}
\leq1/\log n$ and such that $\sum_{i=1}^{k_2} \alpha_i^{(r)}w_{ni}
\geq\beta k_2^{2r/d}$;
\item$\sum_{i=k_2+1}^n w_{ni}^2/\sum_{i=1}^n w_{ni}^2 \leq1/\log n$;
\item$\sum_{i=1}^n |w_{ni}|^3/(\sum_{i=1}^n w_{ni}^2)^{3/2} \leq
1/\log n$.
\end{itemize}
Finally, we are in a position to state the analogue of Theorem~\ref
{ThmMain} for weight vectors in $W_{n,\beta,r}^\dagger$.
%
\begin{thmm}
\label{ThmMain2}
Assume \textup{({A.1})}, \textup{({A.2})({r})}, \textup{({A.3})}
and~\textup{({A.4})}. Then for each $\beta\in(0,1/2)$,
%
\begin{equation}
\label{EqHigherOrderExp} R_{\mathcal{R}}\bigl(\hat{C}_n^{\mathrm{wnn}}
\bigr) - R_{\mathcal{R}}\bigl(C^{\mathrm
{Bayes}}\bigr) = \gamma_n^{(r)}(
\mathbf{w}_n)\bigl\{1 + o(1)\bigr\}
\end{equation}
as $n \rightarrow\infty$, uniformly for $\mathbf{w}_n \in W_{n,\beta
,r}^\dagger$, where
%
\begin{equation}
\label{Eqgamma} \gamma_n^{(r)}(\mathbf{w}_n) =
B_1 \sum_{i=1}^n
w_{ni}^2 + B_2^{(r)} \Biggl(\sum
_{i=1}^n \frac{\alpha_i^{(r)} w_{ni}}{n^{2r/d}}
\Biggr)^2.
\end{equation}
\end{thmm}
A consequence of Theorem~\ref{ThmMain2} is that we can construct
weighted nearest neighbour classifiers which, under conditions
({A.1}), ({A.2})({r}), ({A.3}) and~({A.4}),
and provided that $B_2^{(r)} > 0$, achieve the rate of convergence
$O(n^{-4r/(4r+d)})$ for the regret. To illustrate this,\vadjust{\goodbreak} set $k^{*(r)} =
\lfloor B^{*(r)}n^{4r/(4r+d)} \rfloor$, and in order to satisfy the
restrictions on the allowable weights, consider weight vectors with
$w_{ni} = 0$ for $i= k^{*(r)}+1,\ldots,n$. Then, by mimicking the proof
of Theorem~\ref{ThmOptWeights} and seeking to minimise~(\ref
{Eqgamma}) subject to the constraints $\sum_{i=1}^{k^{*(r)}} w_{ni} =
1$ and $\sum_{i=1}^{k^{*(r)}} \alpha_i^{(\ell)}w_{ni} = 0$ for $\ell=
1,\ldots,r-1$, we obtain minimising weights of the form
%
\begin{equation}
\label{EqNegWeightVec}\quad\hspace*{11pt} w_{ni}^{*(r)} = \cases{\displaystyle
\frac{1}{k^{*(r)}}\bigl(b_0 + b_1
\alpha_i^{(1)} + \cdots+ b_r
\alpha_i^{(r)}\bigr), &\quad  $\mbox{for $i=1,\ldots,k^{*(r)}$},$
\vspace*{2pt}\cr
0, & \quad $\mbox{for $i=k^{*(r)}+1,\ldots,n$.}$ }\hspace*{-10pt}
\end{equation}
The equations $\sum_{i=1}^n w_{ni} = 1$
and $\sum_{i=1}^n \alpha_i^{(\ell)}w_{ni} = 0$
for $\ell=1,\ldots,r-1$ for weight vectors of
the form~(\ref{EqNegWeightVec}) yield $r$ linear equations in the
$r+1$ unknowns $b_0,b_1,\ldots,b_r$. Although these equations can be
solved directly in terms of $b_0$ say, simpler expressions are obtained
by solving asymptotic approximations to these equations. In particular,
since it is an elementary fact that for nonnegative integers~$\ell_1$
and~$\ell_2$,
\[
\sum_{i=1}^k \alpha_i^{(\ell_1)}
\alpha_i^{(\ell_2)} = \frac{(d+2\ell_1)(d+2\ell_2)}{d(d+2\ell_1+2\ell_2)}k^{1+2(\ell_1+\ell_2)/d}\bigl\{
1+O\bigl(k^{-2}\bigr)\bigr\}
\]
as $k \rightarrow\infty$, we can just deal with the dominant terms. As
examples, when $r=1$, we find
\[
b_1 = \frac{1}{(k^{*(1)})^{2/d}}(1-b_0),
\]
and when $r=2$, we should take
\[
b_1 = \frac{1}{(k^{*(2)})^{2/d}} \biggl\{\frac{(d+4)^2}{4} -
\frac
{2(d+4)}{d+2}b_0 \biggr\} \quad\mbox{and}\quad b_2 =
\frac{1 - b_0 -
(k^{*(2)})^{2/d}b_1}{(k^{*(2)})^{4/d}}.
\]
Under the conditions of Theorem~\ref{ThmMain2}, and provided
$B_2^{(r)} > 0$, these weighted nearest neighbour classifiers achieve
the $O(n^{-4r/(4r+d)})$ convergence rate. The choice of $b_0$ involves
a trade-off between the desire to keep the remaining squared bias term
$B_2^{(r)} (\sum_{i=1}^{k^{*(r)}} \frac{\alpha_i^{(r)}
w_{ni}^{*(r)}}{n^{2r/d}} )^2$ small, and the need for it to be
large enough to remain the dominant bias term. This reflects the fact
that the asymptotic results of this section should be applied with some
caution. Besides the discomfort many practitioners might feel in using
negative weights, one would anticipate that rather large sample sizes
would be needed for the leading terms in the asymptotic expansion~(\ref
{EqHigherOrderExp}) to dominate the error terms. This is also the
reason why we do not pursue here methods such as Lepski's method
[\citet
{Lepskii1991}] that adapt to an unknown smoothness level around $\mathcal{S}$.

\section{Empirical performance study}
\label{SecEmpirical}

In this section, we assess the relative empirical performance of the
$k$-nearest neighbour classifier, the optimally weighted nearest
neighbour classifier of Section~\ref{SecMain} and the bagged nearest\vadjust{\goodbreak}
neighbour classifier of Section~\ref{SecBagged} on simulated and real
data sets. We consider four general simulation settings, designed to
exhibit different distributional characteristics:

\textit{Setting} 1: $f_1$ is the density of $d$ independent components,
each having a standard Laplace distribution, and $f_2$ is the density
of the $N_d(\theta,I)$ distribution, where $\theta$ denotes a
$d$-vector of ones.

\textit{Setting} 2: $f_1$ is the density of $d$ independent components,
each having the mixture of normals distribution $\frac{1}{2}N(0,1) +
\frac{1}{2}N(3,2)$. Likewise, $f_2$ is the density of $d$ independent
components, each having a $\frac{1}{2}N(1.5,1) + \frac{1}{2}N(4.5,2)$
distribution.

\textit{Setting} 3: For $d \geq2$, let $\Sigma$ denote the $d \times
d$ Toeplitz matrix whose $j$th entry of its first row is $0.6^{j-1}$.
Set $f_1$ to be the density of the $\frac{1}{2}N_d(0,\Sigma) + \frac
{1}{2}N_d(3\theta,2\Sigma)$ distribution, and $f_2$ to be the density
of the $\frac{1}{2}N_d(3\theta/2,\Sigma) + \frac{1}{2}N_d(9\theta
/2,2\Sigma)$ distribution.

\textit{Setting} 4: Both $f_1$ and $f_2$ are densities of independent
components. For $f_1$, each component has a standard Cauchy density.
For $f_2$, the first $\lfloor d/2 \rfloor$ components also have a
standard Cauchy density, while the last $d - \lfloor d/2 \rfloor$
components have a standard Laplace density.

Setting 1 is a relatively benign classification problem. Setting 2
explores the effect of bimodality, and setting 3 combines bimodal
marginals with dependence between the components. Setting 4 combines
heavy-tailed distributions, a lack of location difference and
introduces components which are irrelevant for classification as
nuisance variables. For each setting, we examined the three sample
sizes $n \in\{50,200,1000\}$, five dimensions $d \in\{1,2,3,5,10\}$
(except for setting 3, where the $d=1$ case was omitted as it is
covered in setting 1) and two prior probabilities $\pi\in\{1/2,2/3\}
$. Thus there were 114 simulation scenarios in total, and we used the
Euclidean norm for computing distances throughout.

In each scenario, we took $\mathcal{R} = \mathbb{R}^d$ and computed the
Bayes risk by Monte Carlo integration. For each data set of size $n$
drawn from the relevant populations, we used a slight variant of a
$5$-fold cross validation algorithm to compute~$\hat{k}$, the number of
neighbours used by the $k$-nearest neighbour classifier. Specifically, we
assigned each observation independently and uniformly at random to one
of five groups, and found the minimiser, denoted $\tilde{k}$, of the
cross-validation risk over a grid of 21 equally spaced points (up to
integer rounding) from 5 to $n/2$. The variant arises from the
observation that this minimiser targets the optimal value of $k$ for a
data set of size $4n/5$. Bearing in mind the expression for the optimal
$k^*$ in~(\ref{Eqkstar}), we therefore set $\hat{k} =  (\frac
{5}{4} )^{4/(d+4)}\tilde{k}$ as an appropriate choice for a data
set of size $n$. The number of positive weights for the optimally
weighted classifier was then chosen to be $\mu(\hat{k})$; cf~(\ref
{Eqmuk}). For the bagged nearest neighbour classifier, we used the
``geometric'' weights given in~(\ref{Eqbagged3}), with $q$ given by
$\hat{q}$ in~(\ref{Eqhatq}). For each data set, we computed the
proportion of misclassifications of $n_{\mathrm{test}} = 1000$
independent test points drawn from the appropriate distribution, and
each simulation was repeated 1000 times to yield estimates of the risks
of each of the three classifiers.

It is computationally convenient to evaluate the distance matrix
between all $n_+ = n+n_{\mathrm{test}}$ points at the outset (even
though some distances will not be used), and this takes $O(n_+^2 d)$
operations when $\|\cdot\|$ is an $\ell_p$-norm. It then takes a
further $O(n_+ n \log n)$ operations to choose $\hat{k}$ and classify
the test points. In particular, the computational requirements are of
the same order of magnitude for both the unweighted and weighted
nearest neighbour classifiers.

An alternative to using a cross-validation method for choosing $\hat
{k}$, as pointed out by an anonymous referee, is to estimate the
constants $B_1$ and $B_2$ in~(\ref{Eqkstar}) directly using a plug-in
approach. We discuss this approach in the supplementary material
[\citet
{Samworth2012}] following the proof of Theorem~\ref{ThmMain2}, but
conclude that it seems awkward to propose a satisfactory algorithm for
estimating $B_1$ and $B_2$ directly, and do not pursue it further here.

Our simulation results are presented in Tables~\ref{TableSimResults1}
and~\ref{TableSimResults2}. To save space, we have omitted the results
for $\pi= 2/3$, which were qualitatively similar. As well as the risks
for the three classifiers, we present in the final two columns
estimates of the regret ratios
%
\begin{equation}
\label{EqDoubleRegretRatio} \frac{R(\hat{C}_{n,\mathbf{w}_n^{\mu(\hat{k})}}^{\mathrm{wnn}}) -
R(C^{\mathrm{Bayes}})}{R(\hat{C}_{n,\hat{k}}^{\mathrm{nn}}) -
R(C^{\mathrm{Bayes}})} \quad\mbox{and}\quad \frac{R(\hat{C}_{n,\hat{q}}^{\mathrm
{bnn}}) - R(C^{\mathrm{Bayes}})}{R(\hat{C}_{n,\hat{k}}^{\mathrm{nn}}) -
R(C^{\mathrm{Bayes}})},
\end{equation}
respectively. Standard errors for these estimates are also given, and
were obtained using the delta method.

\begin{table}
\caption{The estimated risks (multiplied by 100) of the Bayes,
$k$-nearest neighbour, optimally~weighted
nearest neighbour and bagged nearest neighbour classifiers in~settings~1~and~2.~The~final two columns
give the regret ratios defined in~(\protect\ref{EqDoubleRegretRatio}).
Standard errors are given in small script}
\label{TableSimResults1}
\begin{tabular*}{\textwidth}{@{\extracolsep{\fill}}ld{2.3}cccccc@{}}
\hline
$\bolds{d}$ & \multicolumn{1}{c}{\textbf{Bayes}} & $\bolds{n}$ & \phantom{0.}\textbf{knn risk} & \textbf{ownn risk} & \textbf{bnn risk} & \textbf{ownn rr}
& \textbf{bnn rr} \\
\hline
\multicolumn{8}{@{}l}{Setting 1} \\
\phantom{0}1&30.02 & \phantom{00}50 & 33.93$_{0.14}$\phantom{0} & 33.77$_{0.14}$\phantom{0} & 34.71$_{0.14}$\phantom{0} &
0.96$_{0.050}$ & 1.20$_{0.057}$ \\
 &  & \phantom{0}200 & 31.53$_{0.066}$& 31.47$_{0.067}$& 31.72$_{0.075}$& 0.96$_{0.061}$ & 1.10$_{0.071}$ \\
& &1000 & 30.72$_{0.046}$& 30.70$_{0.046}$& 30.72$_{0.046}$&
0.97$_{0.093}$ & 1.00$_{0.094}$ \\ [3pt]
\phantom{0}2&24.21 & \phantom{00}50 & 28.77$_{0.11}$\phantom{0} & 28.58$_{0.11}$\phantom{0} & 28.92$_{0.11}$\phantom{0} &
0.96$_{0.034}$ & 1.00$_{0.035}$ \\
 &  &\phantom{0} 200 & 26.50$_{0.055}$& 26.42$_{0.056}$& 26.51$_{0.058}$& 0.97$_{0.034}$ & 1.00$_{0.035}$ \\
& & 1000& 25.67$_{0.046}$& 25.62$_{0.046}$& 25.60$_{0.046}$&
0.96$_{0.044}$ & 0.95$_{0.044}$ \\ [3pt]
\phantom{0}3&19.37 & \phantom{00}50 & 25.35$_{0.10}$\phantom{0} & 25.03$_{0.098}$& 25.23$_{0.097}$&
0.95$_{0.023}$ & 0.98$_{0.023}$ \\
 &  & \phantom{0}200 & 22.82$_{0.052}$& 22.69$_{0.053}$& 22.72$_{0.053}$& 0.96$_{0.021}$ & 0.97$_{0.021}$ \\
& &1000 & 21.54$_{0.045}$& 21.44$_{0.045}$& 21.43$_{0.046}$&
0.95$_{0.029}$ & 0.95$_{0.029}$ \\ [3pt]
\phantom{0}5& 13.17& \phantom{00}50 & 20.37$_{0.093}$& 20.26$_{0.095}$& 20.40$_{0.093}$&
0.98$_{0.018}$ & 1.00$_{0.018}$ \\
 &  & \phantom{0}200 & 17.74$_{0.049}$& 17.54$_{0.050}$& 17.55$_{0.050}$& 0.96$_{0.015}$ & 0.96$_{0.015}$ \\
& &1000 & 16.21$_{0.041}$& 16.03$_{0.042}$& 16.07$_{0.044}$&
0.94$_{0.019}$ & 0.96$_{0.019}$ \\ [3pt]
10& 5.592& \phantom{00}50 & 13.89$_{0.10}$\phantom{0} & 14.59$_{0.12}$\phantom{0} & 14.63$_{0.11}$\phantom{0} &
1.10$_{0.019}$ & 1.10$_{0.019}$ \\
 &  & \phantom{0}200 & 11.35$_{0.050}$& 11.70$_{0.051}$& 11.72$_{0.050}$& 1.10$_{0.013}$ & 1.10$_{0.013}$ \\
& & 1000& \phantom{00}9.977$_{0.033}$& \phantom{00}9.796$_{0.033}$& \phantom{00}9.911$_{0.033}$&
0.96$_{0.010}$ & 0.98$_{0.011}$ \\[6pt]
\multicolumn{8}{@{}l}{Setting 2} \\
\phantom{0}1& 34.85& \phantom{00}50 & 38.96$_{0.14}$\phantom{0} & 38.78$_{0.13}$\phantom{0} & 39.01$_{0.11}$\phantom{0} &
0.96$_{0.046}$ & 1.00$_{0.043}$ \\
 &  & \phantom{0}200 & 36.76$_{0.073}$& 36.63$_{0.073}$& 36.83$_{0.075}$& 0.93$_{0.052}$ & 1.00$_{0.056}$ \\
& & 1000& 35.34$_{0.052}$& 35.30$_{0.052}$& 35.35$_{0.052}$&
0.91$_{0.14}$\phantom{0} & 1.00$_{0.15}$\phantom{0} \\ [3pt]
\phantom{0}2& 26.83& \phantom{00}50 & 34.36$_{0.13}$\phantom{0} & 33.53$_{0.13}$\phantom{0} & 33.43$_{0.12}$\phantom{0} &
0.89$_{0.023}$ & 0.88$_{0.022}$ \\
 &  & \phantom{0}200 & 30.00$_{0.070}$& 29.63$_{0.068}$& 29.66$_{0.070}$& 0.88$_{0.029}$ & 0.89$_{0.030}$ \\
& &1000 & 27.56$_{0.050}$& 27.48$_{0.050}$& 27.49$_{0.050}$&
0.89$_{0.091}$ & 0.90$_{0.091}$ \\ [3pt]
\phantom{0}3& 21.73& \phantom{00}50 & 31.07$_{0.11}$\phantom{0} & 30.07$_{0.11}$\phantom{0} & 30.01$_{0.11}$\phantom{0} &
0.89$_{0.016}$ & 0.89$_{0.016}$ \\
 &  & \phantom{0}200 & 26.44$_{0.063}$& 25.99$_{0.063}$& 25.96$_{0.065}$& 0.90$_{0.018}$ & 0.90$_{0.018}$ \\
& &1000 & 23.19$_{0.045}$& 23.04$_{0.046}$& 23.03$_{0.046}$&
0.90$_{0.042}$ & 0.90$_{0.042}$ \\ [3pt]
\phantom{0}5& 15.23&\phantom{00}50 & 25.72$_{0.12}$\phantom{0} & 24.88$_{0.11}$\phantom{0} & 25.12$_{0.11}$\phantom{0} &
0.92$_{0.015}$ & 0.94$_{0.015}$ \\
 &  &\phantom{0}200 & 21.51$_{0.055}$& 20.92$_{0.054}$& 20.93$_{0.055}$& 0.91$_{0.012}$ & 0.91$_{0.012}$ \\
& &1000 & 18.69$_{0.045}$& 18.34$_{0.046}$& 18.33$_{0.047}$&
0.90$_{0.018}$ & 0.90$_{0.018}$ \\ [3pt]
10& 7.146&\phantom{00}50 & 16.87$_{0.099}$& 16.57$_{0.10}$\phantom{0} & 16.88$_{0.10}$\phantom{0} &
0.97$_{0.014}$ & 1.00$_{0.015}$ \\
 &  &\phantom{0}200 & 13.00$_{0.048}$& 12.77$_{0.051}$& 12.87$_{0.051}$& 0.96$_{0.012}$ & 0.98$_{0.012}$ \\
& &1000 & 11.57$_{0.034}$& 11.41$_{0.033}$& 11.44$_{0.034}$&
0.96$_{0.011}$ & 0.97$_{0.011}$ \\
\hline
\end{tabular*}   \vspace*{-3pt}
\end{table}

\begin{table}
\caption{The estimated risks (multiplied by 100) of the Bayes,
$k$-nearest neighbour, optimally weighted
nearest neighbour and bagged nearest neighbour classifiers in settings~3
and 4. The final two columns give the regret ratios
defined in~(\protect\ref{EqDoubleRegretRatio}). Standard errors are
given in small script}
\label{TableSimResults2}
\begin{tabular*}{\textwidth}{@{\extracolsep{\fill}}ld{2.3}cccccc@{}}
\hline
$\bolds{d}$ & \multicolumn{1}{c}{\textbf{Bayes}} & $\bolds{n}$ & \phantom{0.}\textbf{knn risk} & \textbf{ownn risk} & \textbf{bnn risk} & \textbf{ownn rr}
& \textbf{bnn rr} \\
\hline
\multicolumn{8}{@{}l}{Setting 3} \\
\phantom{0}2& 32.45& \phantom{00}50 & 37.91$_{0.13}$\phantom{0} & 37.40$_{0.12}$\phantom{0} & 37.54$_{0.11}$\phantom{0} &
0.91$_{0.030}$\phantom{0} & 0.93$_{0.030}$\phantom{0} \\
 &  & \phantom{0}200 & 35.08$_{0.065}$& 34.96$_{0.065}$& 35.05$_{0.068}$& 0.95$_{0.034}$\phantom{0} & 0.99$_{0.036}$\phantom{0} \\
& &1000 & 33.70$_{0.052}$& 33.65$_{0.051}$& 33.67$_{0.051}$&
0.96$_{0.057}$\phantom{0} & 0.98$_{0.057}$\phantom{0} \\ [3pt]
\phantom{0}3& 30.00& \phantom{00}50 & 36.56$_{0.13}$\phantom{0} & 35.94$_{0.11}$\phantom{0} & 36.00$_{0.11}$\phantom{0} &
0.91$_{0.024}$\phantom{0} & 0.91$_{0.024}$\phantom{0} \\
 &  & \phantom{0}200 & 33.61$_{0.065}$& 33.52$_{0.066}$& 33.58$_{0.067}$& 0.97$_{0.025}$\phantom{0} & 0.99$_{0.026}$\phantom{0} \\
& &1000 & 32.03$_{0.052}$& 31.96$_{0.052}$& 31.94$_{0.052}$&
0.96$_{0.036}$\phantom{0} & 0.96$_{0.036}$\phantom{0} \\ [3pt]
\phantom{0}5& 26.13&\phantom{00}50 & 34.10$_{0.14}$\phantom{0} & 33.41$_{0.12}$\phantom{0} & 33.47$_{0.11}$\phantom{0} &
0.91$_{0.021}$\phantom{0} & 0.92$_{0.021}$\phantom{0} \\
 &  &\phantom{0}200 & 30.27$_{0.068}$& 30.16$_{0.070}$& 30.26$_{0.070}$& 0.97$_{0.023}$\phantom{0} & 1.00$_{0.024}$\phantom{0} \\
& &1000 & 28.41$_{0.051}$& 28.23$_{0.051}$& 28.25$_{0.052}$&
0.92$_{0.030}$\phantom{0} & 0.93$_{0.031}$\phantom{0} \\ [3pt]
10& 18.26&\phantom{00}50 & 27.03$_{0.13}$\phantom{0} & 26.50$_{0.11}$\phantom{0} & 26.59$_{0.11}$\phantom{0} &
0.94$_{0.019}$\phantom{0} & 0.95$_{0.019}$\phantom{0} \\
 &  &\phantom{0}200 & 22.86$_{0.067}$& 22.90$_{0.070}$& 23.01$_{0.071}$& 1.00$_{0.021}$\phantom{0} & 1.00$_{0.022}$\phantom{0} \\
& &1000 & 21.07$_{0.046}$& 20.91$_{0.046}$& 20.92$_{0.046}$&
0.94$_{0.022}$\phantom{0} & 0.95$_{0.023}$\phantom{0} \\ [6pt]
\multicolumn{8}{@{}l}{Setting 4} \\
\phantom{0}1& 41.95& \phantom{00}50 & 47.73$_{0.099}$ & 47.49$_{0.10}$\phantom{0}& 47.08$_{0.094}$ &
0.96$_{0.024}$\phantom{0} & 0.89$_{0.022}$\phantom{0} \\
 &  &\phantom{0} 200 & 45.64$_{0.078}$ & 45.45$_{0.077}$& 45.24$_{0.072}$ & 0.95$_{0.029}$\phantom{0} & 0.89$_{0.027}$\phantom{0} \\
& &1000 & 43.38$_{0.061}$ & 43.28$_{0.060}$& 43.32$_{0.061}$
& 0.93$_{0.058}$\phantom{0} & 0.96$_{0.059}$\phantom{0} \\ [3pt]
\phantom{0}2& 41.96& \phantom{00}50 & 48.36$_{0.079}$ & 48.05$_{0.083}$& 47.85$_{0.081}$ &
0.95$_{0.017}$\phantom{0} & 0.92$_{0.017}$\phantom{0} \\
 &  & \phantom{0}200 & 46.39$_{0.074}$ & 46.05$_{0.072}$& 45.96$_{0.070}$ & 0.92$_{0.022}$\phantom{0} & 0.90$_{0.022}$\phantom{0} \\
& &1000 & 44.13$_{0.060}$ & 43.91$_{0.060}$& 43.86$_{0.060}$
& 0.90$_{0.037}$\phantom{0} & 0.88$_{0.037}$\phantom{0} \\ [3pt]
\phantom{0}3& 36.37&\phantom{00}50 & 46.32$_{0.10}$\phantom{0} & 45.73$_{0.10}$\phantom{0} & 45.50$_{0.10}$\phantom{0} &
0.94$_{0.014}$\phantom{0} & 0.92$_{0.014}$\phantom{0} \\
 &  &\phantom{0}200 & 42.92$_{0.083}$ & 42.38$_{0.081}$& 42.29$_{0.078}$ & 0.92$_{0.017}$\phantom{0} & 0.90$_{0.017}$\phantom{0} \\
& &1000 & 39.36$_{0.058}$ & 39.04$_{0.057}$& 39.03$_{0.058}$
& 0.89$_{0.026}$\phantom{0} & 0.89$_{0.026}$\phantom{0} \\ [3pt]
\phantom{0}5& 32.00&50 & 45.66$_{0.10}$\phantom{0} & 44.80$_{0.11}$\phantom{0} & 44.57$_{0.11}$\phantom{0} &
0.94$_{0.011}$\phantom{0} & 0.92$_{0.010}$\phantom{0} \\
 &  &\phantom{0}200 & 40.89$_{0.085}$ & 40.23$_{0.080}$& 40.22$_{0.078}$ & 0.93$_{0.013}$\phantom{0} & 0.93$_{0.012}$\phantom{0} \\
& &1000 & 36.90$_{0.056}$ & 36.45$_{0.056}$& 36.44$_{0.056}$
& 0.91$_{0.015}$\phantom{0} & 0.91$_{0.015}$\phantom{0} \\ [3pt]
10& 25.40&\phantom{00}50 & 45.27$_{0.099}$ & 44.21$_{0.10}$\phantom{0} & 43.97$_{0.098}$ &
0.95$_{0.0069}$& 0.93$_{0.0068}$ \\
 &  &\phantom{0}200 & 39.51$_{0.078}$ & 38.84$_{0.073}$& 38.83$_{0.073}$ & 0.95$_{0.0074}$& 0.95$_{0.0074}$ \\
& &1000 & 36.03$_{0.053}$ & 35.61$_{0.053}$& 35.76$_{0.054}$
& 0.96$_{0.0069}$& 0.97$_{0.0070}$ \\
\hline
\end{tabular*}
\end{table}

In 54 of the 57 scenarios in Tables~\ref{TableSimResults1} and~\ref
{TableSimResults2}, the risk of the optimally weighted nearest
neighbour classifier is smaller than that of the $k$-nearest neighbour
classifier. In one of the three exceptional cases, the difference is so
small that it can easily be explained by the Monte Carlo error. The
other cases are in setting 1 with $d=10$ and $n=50$, $200$. Here it seems
that in this relatively large dimension for nonparametric inference,
these sample sizes are not large enough for the asymptotics to provide
a good approximation.

The extent of the improvement of the optimally weighted nearest
neighbour classifier is generally in close agreement with that predicted
by the theory of Corollary~\ref{CorRegretRatio} and the paragraph
which follows it, even for small sample sizes. This theory tells us
that the first regret ratio in~(\ref{EqDoubleRegretRatio}) converges
to 0.943, 0.924, 0.919, 0.920 and 0.936 in dimensions $d=1,2,3,5,10$,
respectively. Note that a few of the regret ratio estimates,
particularly in settings~1 and~2 with small $d$ and large $n$, have
larger standard errors. This is caused by the fact that in these
scenarios, the risks of all three classifiers are very close to the
Bayes risk. In the more complex situations, the risks of the empirical
classifiers are further from the Bayes risk, and the regret ratios can
be estimated more precisely. The situation is similar for the bagged
nearest neighbour classifier, whose relative performance also matches
that predicted by the theory of Section~\ref{SecBagged} quite
well.\vadjust{\goodbreak}

We also applied all three classifiers to three benchmark data sets,
referred to below as Glass, Yeast and Segmentation, from the UCI
repository [\citet{FrankAsuncion2010}]. Detailed information on these
data sets can be obtained from \url{http://archive.ics.uci.edu/ml/datasets.html}, but summary information
is provided in Table~\ref{TableRealDataResults}. Following \citet
{AthitsosSclaroff2005}, in each case we scaled each component of the
covariates to have unit Euclidean length, and explored both the $\ell_1$ and $\ell_2$ norms for computing distances between observations.
For the Glass and Yeast data sets, we randomly assigned each
observation to a training or test set, each with probability $1/2$,
while for the Segmentation data set, these probabilities were $1/11$
and $10/11$, respectively, since the original data were divided into a
training and test set with these proportions. We then applied the same
modified cross-validation algorithm as for the simulated data to choose
the tuning parameters of the respective procedures. To estimate the
risks of the three classifiers, we computed the proportion of
misclassifications on the test set, and averaged these proportions over
1000 repetitions of the random assignment process.

The results are given in Table~\ref{TableRealDataResults}. In all
cases, the optimally weighted nearest neighbour classifier outperforms
the $k$-nearest neighbour classifier. Since the dimensions for the three
data sets are $d = 9, 8$ and $19$, it is not a surprise to see that the
bagged nearest neighbour classifier also performs comparably well. The
choice of distance appears to make little difference to the relative
performance of the classifiers.

\begin{table}
\caption{The estimated risks (multiplied by 100) of the Bayes,
$k$-nearest neighbour, optimally weighted nearest
neighbour and bagged nearest neighbour classifiers on three UCI
repository data sets. Standard errors are given in small script.
Recall here that $K$ is the number of categories for the response $Y$}
\label{TableRealDataResults}
\begin{tabular*}{\textwidth}{@{\extracolsep{\fill}}lccccccc@{}}
\hline
\textbf{Data set} & \textbf{Distance} & $\bolds{n}$ & $\bolds{d}$ & $\bolds{K}$ & \textbf{knn risk}& \textbf{ownn risk}& \multicolumn{1}{c@{}}{\textbf{bnn risk}}
\\
\hline
Glass & $L_1$ & \phantom{0}163 & \phantom{0}9 & 2 & 23.26$_{0.15}$\phantom{0} & 20.87$_{0.15}$\phantom{0} &
20.36$_{0.15}$\phantom{0} \\
Glass & $L_2$ & \phantom{0}163 & \phantom{0}9 & 2 & 26.21$_{0.15}$\phantom{0} & 23.43$_{0.14}$\phantom{0} &
23.05$_{0.14}$\phantom{0} \\
Yeast & $L_1$ & 1136& \phantom{0}8 & 3 & 40.66$_{0.059}$& 39.71$_{0.062}$&
39.78$_{0.063}$ \\
Yeast & $L_2$ & 1136& \phantom{0}8 & 3 & 40.91$_{0.057}$& 39.90$_{0.058}$&
39.99$_{0.059}$ \\
Segmentation & $L_1$ & 2310 & 19 & 7 & 12.04$_{0.051}$ & 10.05$_{0.043}$ & \phantom{00}9.882$_{0.041}$ \\
Segmentation & $L_2$ & 2310 & 19 & 7 & 15.80$_{0.062}$ & 12.92$_{0.049}$ & 12.67$_{0.049}$ \\
\hline
\end{tabular*}
\end{table}

\begin{appendix}\label{app}
\section*{Appendix}

\begin{pf*}{Proof of Theorem~\ref{ThmMain}}
The proof is rather lengthy, so we briefly outline the main ideas here.
Write $P^\circ= \pi P_1 - (1-\pi)P_2$ and observe that
%
\setcounter{equation}{0}
\begin{eqnarray}
\label{EqMainArg}
&&R_{\mathcal{R}}\bigl(\hat{C}_n^{\mathrm{wnn}}
\bigr) - R_{\mathcal{R}}\bigl(C^{\mathrm
{Bayes}}\bigr)
\nonumber
\\
&&\qquad= \int_{\mathcal{R}} \pi \bigl[\mathbb{P}\bigl\{
\hat{C}_n^{\mathrm{wnn}}(x) = 2\bigr\} - \mathbh{1}_{\{C^{\mathrm{Bayes}}(x)=2\}}
\bigr] \,dP_1(x)
\\
&&\qquad\quad{}+ \int_{\mathcal{R}} (1-\pi) \bigl[\mathbb{P}\bigl\{\hat
{C}_n^{\mathrm{wnn}}(x) = 1\bigr\} - \mathbh{1}_{\{C^{\mathrm{Bayes}}(x)=1\}
}
\bigr] \,dP_2(x)
\nonumber
\\
&&\qquad=\int_{\mathcal{R}} \Biggl\{\mathbb{P} \Biggl(\sum
_{i=1}^n w_{ni} \mathbh{1}_{\{Y_{(i)}=1\}}
< \frac{1}{2} \Biggr) - \mathbh{1}_{\{\eta(x)
< 1/2\}} \Biggr\}
\,dP^\circ(x).\nonumber
\end{eqnarray}
For $\varepsilon> 0$, let
%
\begin{equation}
\label{EqSepseps} \mathcal{S}^{\varepsilon\varepsilon} = \bigl\{x \in\mathbb{R}^d
\dvtx \eta(x) = 1/2 \mbox{ and } \mathrm{dist}(x,\mathcal{S}) < \varepsilon\bigr\},
\end{equation}
where $\mathrm{dist}(x,\mathcal{S}) = \inf_{x_0 \in\mathcal{S}} \|x -
x_0\|$. Moreover, let
\[
\mathcal{S}^\varepsilon= \biggl\{x_0 + t\frac{\dot{\eta}(x_0)}{\|\dot{\eta
}(x_0)\|}
\dvtx x_0 \in\mathcal{S}^{\varepsilon\varepsilon}, |t| < \varepsilon \biggr\}.
\]
The dominant contribution to the integral in~(\ref{EqMainArg}) comes
from $\mathcal{R} \cap\mathcal{S}^{\varepsilon_n}$, where $\varepsilon_n =
n^{-\beta/4d}$. Since the unit vector $\dot{\eta}(x_0)/\|\dot{\eta
}(x_0)\|$ is orthogonal to the tangent space of $\mathcal{S}$ at $x_0$,
we can decompose the integral over $\mathcal{R} \cap\mathcal
{S}^{\varepsilon_n}$ as an integral along $\mathcal{S}$ and an integral in
the perpendicular direction. We then apply a normal approximation to
the integrand to deduce the result. This normal approximation requires
asymptotic expansions to the mean and variance of the sum of
independent random variables in~(\ref{EqMainArg}), and these are
developed in \textit{step} 1 and \textit{step} 2 below, respectively. In
order to retain the flow of the main argument, we concentrate on the
dominant terms in the first five steps of the argument, simply labelling
the many remainder terms as $R_1, R_2,\ldots.$ The sizes of these
remainder terms are controlled in \textit{step} 6 in the supplementary
material [\citet{Samworth2012}], where we also present an additional side
calculation.

\textit{Step} 1: Let $S_n(x) = \sum_{i=1}^n w_{ni} \mathbh{1}_{\{
Y_{(i)}=1\}}$, let $\mu_n(x) = \mathbb{E}\{S_n(x)\}$, let $\varepsilon_n =
n^{-\beta/4d}$ and write $t_n = n^{-2/d}\sum_{i=1}^n \alpha_i w_{ni}$.
We show that
\[
\sup_{x \in\mathcal{S}^{\varepsilon_n}} \bigl|\mu_n(x) - \eta(x) - a(x)t_n x\bigr| =
o(t_n),
\]
uniformly for $\mathbf{w}_n = (w_{ni})_{i=1}^n \in W_{n,\beta}$, where
$a$ is given in~(\ref{Eqa}).

By a Taylor expansion,
%
\begin{eqnarray}
\label{EqTaylor} \mu_n(x) &=& \sum_{i=1}^n
w_{ni}\mathbb{E}\bigl\{\eta(X_{(i)})\bigr\}\nonumber\\
& =& \eta(x) +
\sum_{i=1}^{k_2} w_{ni}\mathbb{E}
\bigl\{(X_{(i)}-x)^T\dot{\eta}(x)\bigr\}
\\
&&{}+ \frac{1}{2}\sum_{i=1}^{k_2}
w_{ni} \mathbb{E}\bigl\{ (X_{(i)}-x)^T\ddot{
\eta}(x) (X_{(i)}-x)\bigr\} + R_1,\nonumber
\end{eqnarray}
where we show in \textit{step} 6 that
%
\begin{equation}
\label{EqR1} \sup_{x \in\mathcal{S}^{\varepsilon_n}} |R_1| = o(t_n),
\end{equation}
uniformly for $\mathbf{w}_n \in W_{n,\beta}$. Writing $p_t = p_t(x) =
\mathbb{P}(\|X-x\| \leq t)$, we also show in \textit{step} 6 that for
$x \in\mathcal{S}^{\varepsilon_n}$ and $i \leq k_2$, the restriction of
the distribution of $X_{(i)} - x$ to\vadjust{\goodbreak} a sufficiently small ball about
the origin is absolutely continuous with respect of Lebesgue measure,
with Radon--Nikodym derivative given at $u = (u_1,\ldots,u_d)^T$ by
%
\begin{eqnarray}
\label{EqDensity} f_{(i)}(u) &=&n\bar{f}(x+u)\pmatrix{{n-1}
\cr
{i-1}}p_{\|u\|}^{i-1}(1-p_{\|
u\|})^{n-i}
\nonumber
\\[-8pt]
\\[-8pt]
\nonumber
 &=& n
\bar{f}(x+u)p_{\|u\|}^{n-1}(i-1),
\end{eqnarray}
say, where $p_{\|u\|}^n(i-1)$ denotes the probability that a $\mathrm
{Bin}(n-1,p_{\|u\|})$ random variable is equal to $i-1$. Let $\delta_n
= (k_2/n)^{1/2d}$. By examining the argument leading to~(0.7) in the supplementary material [\citet{Samworth2012}],
we see that we can replace $\delta$ there with $\delta_n$, to conclude
that for all $M > 0$,
\[
\sup_{x \in\mathcal{S}^{\varepsilon_n}} \sup_{1 \leq i \leq k_2} \mathbb {E}\bigl\{
\|X_{(i)}-x\|^2\mathbh{1}_{\{\|X_{(i)} - x\| > \delta_n\}}\bigr\} = O
\bigl(n^{-M}\bigr).
\]
It follows that
%
\begin{eqnarray}
\label{EqFirstDeriv} \qquad &&\mathbb{E}\bigl\{(X_{(i)}-x)^T\dot{
\eta}(x)\bigr\}
\nonumber
\\[-8pt]
\\[-8pt]
\nonumber
&&\qquad= \int_{\|u\| \leq\delta_n} \dot{\eta}(x)^T u n
\bigl\{\bar{f}(x+u) - \bar{f}(x)\bigr\}p_{\|u\|}^{n-1}(i-1) \,du
+ O\bigl(n^{-M}\bigr),
\end{eqnarray}
uniformly for $x \in\mathcal{S}^{\varepsilon_n}$ and $1 \leq i \leq k_2$.
Similarly,
\begin{eqnarray}
\label{EqSecondDeriv} &&\mathbb{E}\bigl\{(X_{(i)}-x)^T\ddot{
\eta}(x) (X_{(i)}-x)\bigr\}
\nonumber
\\[-8pt]
\\[-8pt]
\nonumber
&&\qquad = \int_{\|u\| \leq
\delta_n}
u^T \ddot{\eta}(x)u n\bar{f}(x+u)p_{\|u\|}^{n-1}(i-1)
\,du + O\bigl(n^{-M}\bigr),
\end{eqnarray}
uniformly for $x \in\mathcal{S}^{\varepsilon_n}$ and $1 \leq i \leq k_2$.
Let $k_1 = \lceil n^{\beta/4} \rceil$, and let $\Delta w_{ni} = w_{ni}
- w_{n,i+1}$ with $w_{n,n+1} = 0$ (where we introduce the comma here
for clarity). By a Taylor expansion, we have
%
\begin{eqnarray}
\label{EqDelta} &&\sum_{i=1}^{k_2}
w_{ni} \int_{\|u\| \leq\delta_n} \biggl[\dot{\eta
}(x)^T u n\bigl\{\bar{f}(x + u) - \bar{f}(x)\bigr\} \nonumber\\
&&\hspace*{68pt}\qquad{}+
\frac{1}{2}u^T \ddot{\eta }(x)u n\bar{f}(x + u)
\biggr]p_{\|u\|}^{n-1}(i - 1) \,du
\nonumber
\\[-8pt]
\\[-8pt]
\nonumber
&&\qquad= \bigl\{1+o(1)\bigr\} \sum_{i=k_1}^{k_2} n
\Delta w_{ni}\sum_{j=1}^d \int
_{\|
u\| \leq\delta_n} \biggl\{\eta_j(x)u_j^2
\bar{f}_j(x) \\
&&\hspace*{166pt}\qquad{}+ \frac{1}{2}\eta_{jj}(x)u_j^2
\bar{f}(x) \biggr\}q_{\|u\|}^{n-1}(i) \,du,\nonumber
\end{eqnarray}
uniformly for $x \in\mathcal{S}^{\varepsilon_n}$ and $\mathbf{w}_n \in
W_{n,\beta}$, where $q_{\|u\|}^{n-1}(i)$ denotes the probability that a
$\mathrm{Bin}(n-1,p_{\|u\|})$ random variable is less than $i$. Now,
$q_{\|u\|}^{n-1}(i)$ is decreasing in\vadjust{\goodbreak} $\|u\|$ and is close to 1 when $\|
u\|$ is small and close to zero when $\|u\|$ is large. To analyse this
more precisely, note that $p_{\|u\|} = \bar{f}(x)a_d \|u\|^d\{1+O(\|u\|^2)\}$ as $u \rightarrow0$,
uniformly for $x \in\mathcal{S}^{\varepsilon
_n}$, so it is convenient to let $b_n =  (\frac{(n-1)a_d \bar
{f}(x)}{i} )^{1/d}$ and set $v = b_n u$. Then there exists $n_0$
such that for $n \geq n_0$, we have for all $x \in\mathcal{S}^{\varepsilon
_n}$, all $\|v\|^d \in(0, 1 - 2/\log n]$ and all $k_1 \leq i \leq k_2$ that
\[
i - (n-1)p_{\|v\|/b_n} \geq\frac{i}{\log n}.
\]
Thus by Bernstein's inequality [\citet{ShorackWellner1986}, page 440],
for each $M > 0$ and for $n \geq n_0$,
%
\begin{equation}
\label{Eqlower} \sup_{\|v\|^d \in(0, 1 - 2/\log n]} \sup_{k_1 \leq i \leq k_2} \bigl\{1 -
q_{\|v\|/b_n}^{n-1}(i)\bigr\} \leq\exp \biggl(-\frac{k_1}{3\log^2 n}
\biggr) = O\bigl(n^{-M}\bigr).\hspace*{-35pt}
\end{equation}
Similarly, for $n \geq n_0$,
%
\begin{equation}
\label{Equpper}\qquad \sup_{\|v\|^d \in[1 + 2/\log n,b_n\delta_n]} \sup_{k_1 \leq i \leq
k_2} q_{\|v\|/b_n}^{n-1}(i)
\leq\exp \biggl(-\frac{k_1}{3\log^2 n} \biggr) = O\bigl(n^{-M}\bigr).\hspace*{-35pt}
\end{equation}
We deduce from (\ref{EqFirstDeriv})--(\ref{Eqlower}) and (\ref{Equpper}) that
%
\begin{eqnarray}
\label{EqBiascompletion} &&\sum_{i=1}^{k_2}
w_{ni}\mathbb{E}\bigl\{(X_{(i)}-x)^T\dot{
\eta}(x)\bigr\} \nonumber\\
&&\quad{}+ \frac
{1}{2}\sum_{i=1}^{k_2}
w_{ni} \mathbb{E}\bigl\{(X_{(i)}-x)^T\ddot{\eta
}(x) (X_{(i)}-x)\bigr\}
\nonumber\\
&&\qquad= \bigl\{1 + o(1)\bigr\}\\
&&\qquad\quad{}\times\sum_{i=1}^n
\frac{n \Delta w_{ni}}{b_n^{d+2}} \sum_{j=1}^d \biggl\{
\eta_j(x)\bar{f}_j(x) + \frac{1}{2}
\eta_{jj}(x)\bar {f}(x) \biggr\} \int_{\|v\| \leq1}
v_j^2 \,dv
\nonumber\\
&&\qquad= a(x)t_n + o(t_n),\nonumber
\end{eqnarray}
uniformly for $x \in\mathcal{S}^{\varepsilon_n}$ and $\mathbf{w}_n \in
W_{n,\beta}$. Combining~(\ref{EqTaylor}),~(\ref{EqR1}) and~(\ref
{EqBiascompletion}), this completes \textit{step}~1.

\textit{Step} 2: Let $\sigma_n^2(x) = \mathrm{Var}\{S_n(x)\}$ and let
$s_n^2 = \sum_{i=1}^n w_{ni}^2$. We claim that
\[
\sup_{x \in\mathcal{S}^{\varepsilon_n}} \biggl|\sigma_n^2(x) -
\frac
{1}{4}s_n^2 \biggr| = o\bigl(s_n^2
\bigr),
\]
uniformly for $\mathbf{w}_n \in W_{n,\beta}$. To see this, note that
\begin{eqnarray*}
\sigma_n^2(x) &=& \sum_{i=1}^n
w_{ni}^2 \mathbb{E} \bigl[\eta(X_{(i)})\bigl\{1
- \eta(X_{(i)})\bigr\} \bigr] + \sum_{i=1}^n
w_{ni}^2 \mathrm{Var} \eta (X_{(i)})
\\
&=& \sum_{i=1}^n w_{ni}^2
\bigl[\mathbb{E}\eta(X_{(i)}) - \bigl\{\mathbb {E}\eta(X_{(i)})
\bigr\}^2 \bigr].
\end{eqnarray*}
But by a simplified version of the argument in \textit{step} 1, we have
\[
\sup_{x \in\mathcal{S}^{\varepsilon_n}} \sup_{1 \leq i \leq k_2} \bigl|\mathbb {E}\eta(X_{(i)}) -
\eta(x)\bigr| \rightarrow0.
\]
It follows that
\begin{eqnarray*}
&&\sup_{x \in\mathcal{S}^{\varepsilon_n}} \Biggl|\sum_{i=1}^n
w_{ni}^2 \mathbb{E}\eta(X_{(i)}) -
\frac{1}{2}s_n^2\Biggr |\\[-2pt]
&&\qquad \leq\sup_{x \in
\mathcal{S}^{\varepsilon_n}} \sum
_{i=1}^{k_2} w_{ni}^2
\bigl|\mathbb{E}\eta (X_{(i)}) - \eta(x)\bigr| \\[-2pt]
&&\qquad\quad{}+ \sum
_{i=k_2+1}^n w_{ni}^2
+ s_n^2 \sup_{x \in\mathcal{S}^{\varepsilon_n}} \bigl|\eta(x) - 1/2\bigr| \\[-2pt]
&&\qquad= o
\bigl(s_n^2\bigr),
\end{eqnarray*}
uniformly for $\mathbf{w}_n \in W_{n,\beta}$. Similarly,
\begin{eqnarray*}
&&\Biggl|\sum_{i=1}^n w_{ni}^2
\bigl\{\mathbb{E}\eta(X_{(i)})\bigr\}^2 -
\frac
{1}{4}s_n^2 \Biggr|\nonumber\\[-2pt]
&&\qquad \leq\sum
_{i=1}^{k_2} w_{ni}^2 \bigl|
\mathbb{E}\eta (X_{(i)}) - \eta(x)\bigr|\bigl|\mathbb{E}\eta(X_{(i)}) +
\eta(x)\bigr|
\\[-2pt]
&&\qquad\quad{}+ 2\sum_{i=k_2+1}^n w_{ni}^2
+ s_n^2 \bigl|\eta(x)^2 - 1/4\bigr| \\[-2pt]
&&\qquad= o
\bigl(s_n^2\bigr),\nonumber
\end{eqnarray*}
uniformly for $x \in\mathcal{S}^{\varepsilon_n}$ and $\mathbf{w}_n \in
W_{n,\beta}$. This completes \textit{step} 2.

\textit{Step} 3: For $x_0 \in\mathcal{S}$ and $t \in\mathbb{R}$, we
write $x_0^t = x_0 + t\dot{\eta}(x_0)/\|\dot{\eta}(x_0)\|$ for brevity.
Moreover, we write $\psi= \pi f_1 - (1-\pi)f_2$ for the Radon--Nikodym
derivative with respect to Lebesgue measure of the restriction of
$P^\circ$ to $\mathcal{S}^{\varepsilon_n}$ for large~$n$. We show that
%
\begin{eqnarray}
\label{EqRFubini} &&\int_{\mathcal{R} \cap\mathcal{S}^{\varepsilon_n}} \bigl[\mathbb{P}\bigl\{
S_n(x) < 1/2\bigr\} - \mathbh{1}_{\{\eta(x) < 1/2\}} \bigr]
\,dP^\circ(x)
\nonumber\\[-2pt]
&&\qquad= \int_{\mathcal{S}} \int_{-\varepsilon_n}^{\varepsilon_n}
\psi \bigl(x_0^t\bigr) \bigl[\mathbb{P}\bigl
\{S_n\bigl(x_0^t\bigr) < 1/2\bigr\}\\[-2pt]
&&\hspace*{112pt}\qquad{} -
\mathbh{1}_{\{t < 0\}
} \bigr] \,dt \,d\mathrm{Vol}^{d-1}(x_0)
\bigl\{1 + o(1)\bigr\},\nonumber
\end{eqnarray}
uniformly for $\mathbf{w}_n \in W_{n,\beta}$. Recalling the definition
of $\mathcal{S}^{\varepsilon_n \varepsilon_n}$ in~(\ref{EqSepseps}), note
that for large $n$, the map
\[
\phi \biggl(x_0,t\frac{\dot{\eta}(x_0)}{\|\dot{\eta}(x_0)\|} \biggr) = x_0^t\vadjust{\goodbreak}
\]
is a diffeomorphism from $\{(x_0,t\dot{\eta}(x_0)/\|\dot{\eta}(x_0)\|
)\dvtx x_0 \in\mathcal{S}^{\varepsilon_n \varepsilon_n}, |t| < \varepsilon_n\}$ onto
$\mathcal{S}^{\varepsilon_n}$ [\citet{Gray2004}, pages~32--33]. Observe that
%
\begin{equation}
\label{EqContainment} \bigl\{x \in\mathbb{R}^d\dvtx \mathrm{dist}(x,
\mathcal{S}) < \varepsilon_n\bigr\} \subseteq\mathcal{S}^{\varepsilon_n}
\subseteq\bigl\{x \in\mathbb {R}^d\dvtx \mathrm{dist}(x,\mathcal{S})
< 2\varepsilon_n\bigr\}.
\end{equation}
Moreover, for large $n$ and $|t| < \varepsilon_n$, we have $\mathrm{sgn}\{
\eta(x_0^t)-1/2\} = \mathrm{sgn}\{\psi(x_0^t)\} = \mathrm{sgn}(t)$. The
pullback of the $d$-form $dx$ is given at $(x_0,t\dot{\eta}(x_0)/\|\dot
{\eta}(x_0)\|)$ by
\[
\det\dot{\phi} \biggl(x_0,t\frac{\dot{\eta}(x_0)}{\|\dot{\eta}(x_0)\|
} \biggr) \,dt \,d
\mathrm{Vol}^{d-1}(x_0) = \bigl\{1 + o(1)\bigr\} \,dt \,d
\mathrm{Vol}^{d-1}(x_0),
\]
where the error term is uniform in $(x_0,t\dot{\eta}(x_0)/\|\dot{\eta
}(x_0)\|)$ for $x_0 \in\mathcal{S}$ and $|t| < \varepsilon_n$. It follows
from the theory of integration on manifolds, as described
in~\citet{GuilleminPollack1974},
page~168 and \citet{Gray2004}, Theorems~3.15 and 4.7 [see also \citet{Moore1992}], that
%
\begin{eqnarray}
\label{EqFubini}\qquad &&\int_{\mathcal{S}^{\varepsilon_n}} \bigl[\mathbb{P}\bigl
\{S_n(x) < 1/2\bigr\} - \mathbh{1}_{\{\eta(x) < 1/2\}} \bigr]
\,dP^\circ(x)
\nonumber
\\
&&\qquad= \int_{\mathcal{S}^{\varepsilon_n \varepsilon_n}} \int_{-\varepsilon_n}^{\varepsilon_n}
\psi\bigl(x_0^t\bigr) \bigl[\mathbb{P}\bigl
\{S_n\bigl(x_0^t\bigr) < 1/2\bigr\} \\
&&\hspace*{126pt}\qquad{}-
\mathbh{1}_{\{t < 0\}} \bigr] \,dt \,d\mathrm{Vol}^{d-1}(x_0)
\bigl\{1 + o(1)\bigr\},\nonumber
\end{eqnarray}
uniformly for $\mathbf{w}_n \in W_{n,\beta}$. But $\mathcal{S}^{\varepsilon
_n} \setminus\mathcal{R} \subseteq\{x \in\mathbb{R}^d\dvtx \mathrm
{dist}(x,\partial\mathcal{S}) < \varepsilon_n\}$, and this latter set has
volume $O(\varepsilon_n^2)$ by Weyl's tube formula [\citet{Gray2004}, Theorem~4.8].
Thus the integral over $\mathcal{S}^{\varepsilon
_n}$ in~(\ref{EqFubini}) may be replaced with an integral over
$\mathcal{R} \cap\mathcal{S}^{\varepsilon_n}$ and, similarly, the
integral over $\mathcal{S}^{\varepsilon_n \varepsilon_n}$ may be replaced
with an integral over $\mathcal{S}$, without changing the order of the
error term in~(\ref{EqFubini}). Thus~(\ref{EqRFubini}) holds, and
this completes \textit{step} 3.

\textit{Step} 4: We now return to the main argument to bound the
contribution to the risk~(\ref{EqMainArg}) from $\mathcal{R} \setminus
\mathcal{S}^{\varepsilon_n}$. In particular, we show that
%
\begin{equation}
\label{EqOuter} \sup_{\mathbf{w}_n \in W_{n,\beta}} \int_{\mathcal{R} \setminus\mathcal
{S}^{\varepsilon_n}} \bigl[
\mathbb{P}\bigl\{S_n(x) < 1/2\bigr\} - \mathbh{1}_{\{\eta
(x) < 1/2\}}
\bigr] \,dP^\circ(x) = O\bigl(n^{-M}\bigr)\hspace*{-35pt}
\end{equation}
for all $M > 0$. To see this, recall that $|\eta(x) - 1/2|$ is assumed
to be bounded away from zero on the set $\mathcal{R} \setminus\mathcal
{S}^{\varepsilon}$ (for fixed $\varepsilon> 0$), and $\|\dot{\eta}(x_0)\|$
is bounded away from zero for $x_0 \in\mathcal{S}$. Hence, by~(\ref
{EqContainment}) in \textit{step} 3, there exists $c_1 > 0$ such that,
for sufficiently small $\varepsilon> 0$,
%
\begin{equation}
\label{Eqinfeta} \inf_{x \in\mathcal{R} \setminus\mathcal{S}^{\varepsilon}} \bigl|\eta(x) - 1/2\bigr| \geq c_1
\varepsilon.
\end{equation}
We also claim that $\mu_n(x) = \mathbb{E}\{S_n(x)\}$ is similarly
bounded away from $1/2$ uniformly for $x \in\mathcal{R} \setminus
\mathcal{S}^{\varepsilon_n}$. In fact, we have by Hoeffding's inequality that
\[
\mathbb{P}\bigl(\|X_{(k_2)} - x\| > \varepsilon_n/2\bigr) =
q_{\varepsilon_n/2}^n(k_2) \leq e^{-({2}/{n})(n p_{\varepsilon_n/2} - k_2)^2}.
\]
It follows that
%
\begin{eqnarray}
\label{EqUpper} &&\mathop{\sup_{x \in\mathcal{R} \setminus\mathcal{S}^{\varepsilon
_n}\dvtx }}_{\eta(x) \leq1/2}
\mu_n(x) - \frac{1}{2}\nonumber\\
&&\qquad\leq\mathop{\sup_{x \in\mathcal{R} \setminus\mathcal{S}^{\varepsilon
_n}\dvtx }}_{\eta(x) \leq1/2}
\Biggl\{\sum_{i=1}^{k_2} w_{ni}
\mathbb {P}\bigl(Y_{(i)} = 1 \cap\|X_{(k_2)} - x\| \leq
\varepsilon_n/2\bigr)
\\
&&\hspace*{95pt}{}- \frac{1}{2} + e^{-({2}/{n})(n p_{\varepsilon_n/2} -
k_2)^2} + \sum_{i=k_2+1}^n
w_{ni} \Biggr\}
\nonumber\\
&&\qquad\leq\sum_{i=1}^{k_2} w_{ni}
\biggl(\frac{1}{2} - \frac{c_1\varepsilon_n}{2} \biggr) - \frac{1}{2} +
e^{-({2}/{n})(n p_{\varepsilon_n/2} -
k_2)^2} + n^{-\beta} \leq-\frac{c_1\varepsilon_n}{4}
\end{eqnarray}
for sufficiently large $n$. Similarly,
%
\begin{eqnarray}
\label{EqLower} &&\mathop{\inf_{x \in\mathcal{R} \setminus\mathcal{S}^{\varepsilon
_n}\dvtx }}_{\eta(x) \geq1/2} \mu_n(x) - \frac{1}{2} \nonumber\\
&&\qquad
\geq\mathop{\inf_{x
\in\mathcal{R} \setminus\mathcal{S}^{\varepsilon_n}\dvtx }}_{\eta(x) \geq1/2} \sum_{i=1}^{k_2}
w_{ni}\mathbb{P}\bigl(Y_{(i)} = 1 \cap\|X_{(k_2)} - x\|
\leq\varepsilon_n/2\bigr) - \frac{1}{2}
\\
&&\qquad\geq\bigl(1 - n^{-\beta/2}\bigr) \biggl(\frac{1}{2} +
\frac{c_1\varepsilon_n}{2} \biggr) \bigl(1 - e^{-({2}/{n})(n p_{\varepsilon_n/2} - k_2)^2}\bigr) - \frac{1}{2}
\geq \frac{c_1\varepsilon_n}{4}\nonumber
\end{eqnarray}
for large $n$.

Now we may apply Hoeffding's inequality again, this time to $S_n(x)$,
to deduce that
\[
\bigl|\mathbb{P}\bigl\{S_n(x) < 1/2\bigr\} - \mathbh{1}_{\{\eta(x) < 1/2\}}\bigr |
\leq e^{{(-2(\mu_n(x) - 1/2)^2)}/{s_n^2}} = O\bigl(n^{-M}\bigr)
\]
for each $M > 0$, uniformly for $\mathbf{w}_n \in W_{n,\beta}$ and $x
\in\mathcal{R} \setminus\mathcal{S}^{\varepsilon_n}$, using~(\ref
{EqUpper}) and~(\ref{EqLower}) and the fact that $s_n^2 \leq n^{-\beta
}$ for $\mathbf{w}_n \in W_{n,\beta}$. This completes \textit{step}~4.

\textit{Step} 5. We now show that
\begin{eqnarray*}
&&\int_{\mathcal{S}} \int_{-\varepsilon_n}^{\varepsilon_n}
\psi\bigl(x_0^t\bigr) \bigl[\mathbb{P}\bigl
\{S_n\bigl(x_0^t\bigr) < 1/2\bigr\} -
\mathbh{1}_{\{t < 0\}} \bigr] \,dt\,  d\mathrm{Vol}^{d-1}(x_0)
\\
&&\qquad= B_1 s_n^2 + B_2t_n^2
+ o\bigl(s_n^2 + t_n^2\bigr),
\end{eqnarray*}
uniformly for $\mathbf{w}_n \in W_{n,\beta}$, where $B_1$ and $B_2$
were defined in~(\ref{EqC1C2}). When combined with~(\ref{EqMainArg})
and the results of \textit{step} 3 and \textit{step} 4 [in
particular,~(\ref{EqRFubini}) and~(\ref{EqOuter})], this will
complete the proof of Theorem~\ref{ThmMain}.

First observe that
%
\begin{eqnarray}
\label{EqLinearise} &&\int_{\mathcal{S}} \int_{-\varepsilon_n}^{\varepsilon_n}
\psi\bigl(x_0^t\bigr) \bigl[\mathbb{P}\bigl
\{S_n\bigl(x_0^t\bigr) < 1/2\bigr\} -
\mathbh{1}_{\{t < 0\}} \bigr] \,dt \,d\mathrm{Vol}^{d-1}(x_0)
\nonumber
\\
&&\qquad= \int_{\mathcal{S}} \int_{-\varepsilon_n}^{\varepsilon_n}
t\bigl\|\dot{\psi }(x_0)\bigr\| \bigl[\mathbb{P}\bigl\{S_n
\bigl(x_0^t\bigr) < 1/2\bigr\} \\
&&\hspace*{126pt}\qquad{}- \mathbh{1}_{\{t < 0\}
}
\bigr] \,dt \,d\mathrm{Vol}^{d-1}(x_0)\bigl\{1+o(1)\bigr
\}.\nonumber
\end{eqnarray}
Now, $S_n(x)$ is a sum of independent, bounded random variables, so by
the nonuniform version of the Berry--Esseen theorem, there exists $C_1
> 0$ such that for all $y \in\mathbb{R}$,
\begin{eqnarray*}
&&\sup_{x_0 \in\mathcal{S}} \sup_{t \in[-\varepsilon_n,\varepsilon_n]} \biggl|\mathbb{P} \biggl(\frac{S_n(x_0^t) - \mu_n(x_0^t)}{\sigma_n(x_0^t)}
\leq y \biggr) - \Phi(y) \biggr| \\
&&\qquad\leq\frac{C_1}{n^{1/2}(1+|y|^3)},
\end{eqnarray*}
where $\Phi$ denotes the standard normal distribution function. Thus
\begin{eqnarray*}
&&\int_{\mathcal{S}} \int_{-\varepsilon_n}^{\varepsilon_n} t
\bigl\|\dot{\psi}(x_0)\bigr\| \bigl\{\mathbb{P}\bigl\{S_n
\bigl(x_0^t\bigr) < 1/2\bigr\} - \mathbh{1}_{\{t < 0\}}
\bigr\} \,dt \,d\mathrm{Vol}^{d-1}(x_0)
\\
&&\qquad= \int_{\mathcal{S}} \int_{-\varepsilon_n}^{\varepsilon_n}
t\bigl\|\dot{\psi }(x_0)\bigr\| \biggl\{\Phi \biggl(\frac{1/2 - \mu_n(x_0^t)}{\sigma_n(x_0^t)}
\biggr) \\
&&\hspace*{138pt}\qquad{}- \mathbh{1}_{\{t < 0\}} \biggr\} \,dt \,d\mathrm {Vol}^{d-1}(x_0)
+ R_2,
\end{eqnarray*}
where we show in \textit{step} 6 that
%
\begin{equation}
\label{EqR2} |R_2| = o\bigl(s_n^2 +
t_n^2\bigr),
\end{equation}
uniformly for $\mathbf{w}_n \in W_{n,\beta}$. Moreover, by a Taylor
expansion and \textit{step} 1 and \textit{step} 2,
\begin{eqnarray*}
&&\int_{\mathcal{S}} \int_{-\varepsilon_n}^{\varepsilon_n} t
\bigl\|\dot{\psi}(x_0)\bigr\| \biggl\{\Phi \biggl(\frac{1/2 - \mu_n(x_0^t)}{\sigma_n(x_0^t)} \biggr) -
\mathbh{1}_{\{t < 0\}} \biggr\} \,dt \,d\mathrm{Vol}^{d-1}(x_0)
\\
&&\qquad= \int_{\mathcal{S}} \int_{-\varepsilon_n}^{\varepsilon_n}
t\bigl\|\dot{\psi }(x_0)\bigr\| \biggl\{\Phi \biggl(\frac{-2t\|\dot{\eta}(x_0)\| -
2a(x_0)t_n}{s_n}
\biggr)\\
&&\hspace*{180pt}\qquad{} - \mathbh{1}_{\{t < 0\}} \biggr\} \,dt \,d\mathrm {Vol}^{d-1}(x_0)
+ R_3,
\end{eqnarray*}
where we show in \textit{step} 6 that
%
\begin{equation}
\label{EqR3} |R_3| = o\bigl(s_n^2 +
t_n^2\bigr),
\end{equation}
uniformly for $\mathbf{w}_n \in W_{n,\beta}$. Finally, we can make the
substitution $r = t/s_n$ to conclude that
\begin{eqnarray*}
&&\hspace*{-4pt}\int_{\mathcal{S}} \int_{-\varepsilon_n}^{\varepsilon_n} t
\bigl\|\dot{\psi}(x_0)\bigr\| \biggl\{\Phi \biggl(\frac{-2t\|\dot{\eta}(x_0)\| - 2a(x_0)t_n}{s_n} \biggr) -
\mathbh{1}_{\{t < 0\}} \biggr\} \,dt \,d\mathrm{Vol}^{d-1}(x_0)
\nonumber
\\
&&\hspace*{-4pt}\qquad= \frac{s_n^2}{4} \int_{\mathcal{S}} \int_{-\infty}^\infty
u\bigl\|\dot {\psi}(x_0)\bigr\| \biggl\{\Phi \biggl(-u\bigl\|\dot{
\eta}(x_0)\bigr\| - \frac{2t_n
a(x_0)}{s_n} \biggr)\\
&&\hspace*{189pt}\qquad{} - \mathbh{1}_{\{u < 0\}}
\biggr\} \,du \,d\mathrm {Vol}^{d-1}(x_0)+ R_4
\nonumber
\\
&&\hspace*{-4pt}\qquad= B_1 s_n^2 + B_2t_n^2
+ R_4,
\end{eqnarray*}
where $B_1$ and $B_2$ were defined in~(\ref{EqC1C2}). Here, we have
used the fact that $\|\dot{\psi}(x_0)\|/\|\dot{\eta}(x_0)\| = 2\bar
{f}(x_0)$ for $x_0 \in\mathcal{S}$ in the final step of this
calculation. Once we have shown in \textit{step} 6 that
%
\begin{equation}
\label{EqR4} |R_4| = o\bigl(s_n^2\bigr),
\end{equation}
uniformly for $\mathbf{w}_n \in W_{n,\beta}$, this will complete \textit{Step} 5 and hence the proof of Theorem~\ref{ThmMain}.
\end{pf*}
\end{appendix}

\section*{Acknowledgements}
I am very grateful to Peter Hall for introducing me to this topic, and
to the anonymous
referees and Associate Editor for their constructive comments, which
significantly helped to improve the paper.
I thank the Statistical and Applied Mathematical Sciences Institute
(SAMSI) in North Carolina for kind hospitality when
I attended the ``Program on Analysis of Object Data'' from September to
November 2010,
during which time part of this research was carried out.

\begin{supplement}
\stitle{Supplement to ``Optimal weighted nearest neighbour classifiers''}
\slink[doi]{10.1214/12-AOS1049SUPP} 
\sdatatype{.pdf}
\sfilename{aos1049\_supp.pdf}
\sdescription{We complete the proof of Theorem~\ref{ThmMain}, and give
the proofs of
the other results in the paper. We also discuss minimax properties of
weighted nearest
neighbour classifiers and a plug-in approach to estimating $k^*$.}
\end{supplement}

%

\printaddresses

\end{document}